\documentclass[11pt]{article}

\usepackage{amsmath,amsthm,amssymb,color,epic,eepic,
cite}
\usepackage{amsfonts}

\renewcommand{\baselinestretch}{1.2}

\def\baselinestretch{1.4}
\newlength{\minitwocolumn}
\newcommand{\Z}{{\Bbb Z}} 
\newcommand{\R}{{\Bbb R}} 
\newcommand{\C}{{\Bbb C}} 

\newcommand{\K}{{\Bbb F}} 
\newcommand{\N}{{\Bbb N}} 

\newcommand{\F}{{\cal F}}

\newcommand{\cC}{{\mathfrak C}}
\newcommand{\cD}{{\cal D}}
\newcommand{\fC}{{\mathfrak C}}
\newcommand{\fD}{{\mathfrak D}}
\newcommand{\sC}{{\mathfrak C}}
\newcommand{\sD}{{\mathfrak D}}
\newcommand{\cA}{{\cal A}}
\newcommand{\cH}{{\cal H}}

\newcommand{\cM}{{\cal M}}
\newcommand{\cN}{{\cal N}}

\newcommand{\cZ}{{\cal Z}}

\newcommand{\hcZ}{{\cal Z}}

\newcommand{\cR}{{\cal R}}
\newcommand{\cP}{{\cal P}}
\newcommand{\cQ}{{\cal Q}}
\newcommand{\cE}{{\cal E}}
\newcommand{\cW}{{\cal W}}

\newcommand{\la}{\lambda}

\newcommand{\al}{\alpha}

\newcommand{\ep}{\epsilon}

\newcommand{\bep}{{\epsilon}}
\newcommand{\s}{{\sigma}}

\newcommand{\hd}{{d}}
\newcommand{\hf}{\widehat{f}}
\newcommand{\hs}{{\sigma}}
\newcommand{\bs}{\bar{\sigma}}

\newcommand{\hpi}{{\pi}}

\newcommand{\hV}{{\cal V}}

\newcommand{\bh}{{\bar{\h}}}
\newcommand{\bd}{\bar{d}}

\newcommand{\wt}[1]{{{\rm wt}(#1)}}

\newcommand{\nn}{{\nonumber}}
\newcommand{\bea}{\begin{eqnarray}}
\newcommand{\ena}{\end{eqnarray}}
\newcommand{\beit}{\begin{itemize}}
\newcommand{\enit}{\end{itemize}}

\newcommand{\be}{\begin{eqnarray*}}
\newcommand{\en}{\end{eqnarray*}}
\newcommand{\lb}[1]{\label{#1}}

\newcommand{\ds}[1]{{\displaystyle #1 }}

\newcommand{\End}{{\rm End}}
\newcommand{\Ind}{{\rm Ind}}

\def\infq4p#1{{(#1;q^4,p)_\infty}}

\newcommand{\hPsi}{\widehat{\Psi}}

\newcommand{\tr}{{\rm tr}}

\newcommand{\mmatrix}[1]{\begin{matrix} #1 \end{matrix}}

\font\teneufm=eufm10
\font\seveneufm=eufm7
\font\fiveeufm=eufm5
\newfam\eufmfam
\textfont\eufmfam=\teneufm
\scriptfont\eufmfam=\seveneufm
\scriptscriptfont\eufmfam=\fiveeufm

\let\goth\frak
\newcommand{\slth}{\widehat{\goth{sl}}_2}
\newcommand{\slt}{\goth{sl}_2}
\newcommand{\gsl}{\goth{sl}}
\newcommand{\slnh}{\widehat{\goth{sl}}_N}

\newcommand{\g}{\goth{g}}

\newcommand{\Bqla}{{{\cal B}_{q,\lambda}}}

\newcommand{\h}{\goth{h}}

\newcommand{\gh}{\widehat{\goth{g}}}

\font\fourteeneufm=eufm10 scaled\magstep2    
   
\newcommand{\ghbig}{\widehat{\mbox{\fourteeneufm g}}}  
  
\makeatletter
\@addtoreset{equation}{section}
\makeatother

\newtheorem{thm}{Theorem}[section]
\newtheorem{prop}[thm]{Proposition}

\newtheorem{conj}[thm]{Conjecture}

\newtheorem{dfn}[thm]{Definition}

\begin{document}

\begin{center}
{\Large \bf Elliptic Algebra $U_{q,p}(\ghbig)$ and Quantum $Z$-algebras}\\[10mm]
{\large Rasha M. Farghly$^{*,}$\footnote{E-mail: d102349@hiroshima-u.ac.jp },  Hitoshi Konno$^{*,}$\footnote{Address after April 1, 2014: Dept.of Math. Tokyo University of Marine Science and Technology, Etchujima,\par \quad Tokyo 135-8533, Japan; E-mail: hkonno0@kaiyodai.ac.jp } and Kazuyuki Oshima$^{\star,}$\footnote{ E-mail: oshima@aitech.ac.jp }}\\[5mm]
{\it $^*$Department of Mathematics, Graduate School of Science, Hiroshima University, Higashi-Hiroshima 739-8521, Japan}\\
{\it $^\star$Department of Mathematics, Center for General Education, Aichi Institute of Technology, Yakusa-cho, Toyota
470-0392, Japan}\\

\end{center}

\begin{abstract}
A new definition of the elliptic algebra $U_{q,p}(\gh)$ 
associated with an untwisted affine Lie algebra $\gh$ is given 
as a topological algebra 
over the ring of formal power series in $p$. 
We also introduce a quantum dynamical 
analogue of Lepowsky-Wilson's $Z$-algebras. 
The $Z$-algebra governs the irreducibility of 
 the infinite dimensional $U_{q,p}({\gh})$-modules. 
Some level-1 examples indicate a direct connection of the irreducible $U_{q,p}(\gh)$-modules to those of the 
$W$-algebras associated with the coset 
$\gh\oplus \gh\supset (\gh)_{{diag}}$ with level $(r-g-1,1)$ ($g$:the dual Coxeter number), which includes Fateev-Lukyanov's $W\! B_l$-algebra.  
\end{abstract}

\section{Introduction}

The algebra $U_{q,p}(\gh)$ is an elliptic analogue  
\cite{Konno98,JKOS} of  the quantum affine algebra $U_q(\gh)$ in the Drinfeld realization\cite{Dr}. There are two types of the elliptic quantum groups, the vertex 
type and the face type\cite{Fro,JKOStg}. 
Deriving the $L$-operators\cite{EF, JKOS,KojimaKonno} and 
introducing the Hopf-algebroid structure\cite{EV,KR,Konno09} $U_{q,p}(\gh)$ is now recognized as a face type elliptic quantum group. 

Originally $U_{q,p}(\gh)$ with $p=q^{2r}$ was derived for $\slth=\widehat{\gsl}(2,\C)$\cite{Konno98} 
as a deformation of the screening currents of the coset conformal field theory (CFT) $\slth\oplus  \slth \supset (\slth)_{{ diag}}$ with level $(r-k-2, k)$ \cite{GKO,Rav, KMQ,BNY,DSZ,GMM} instead of   considering a deformation of $U_q(\slth)$ itself. 
 Such coset CFT is known to be realized in terms of the level-$k$ free boson and the $\Z_k$-parafermion\cite{ZamFat}, or the $Z$-algebra\cite{LW} associated with the level-$k$ standard representation of $\slth$\footnote{The difference between $Z$-algebra and Parafermion 
is whether one adds zero-modes 
of the bosons to it or not.}.  
It was then crucial in \cite{Konno98} 
 to realize that the level-$k$ boson should be deformed 
both $q$- and elliptically\cite{LP} whereas the $\Z_k$-parafermion gets 
only a $q$-deformation to obtain {consistent relations} for the generators in $U_{q,p}(\slth)$.

In \cite{JKOS}, a realization of  $U_{q,p}(\gh)$ for general untwisted affine Lie 
algebra $\gh$ was given by modifying the Drinfeld realization of the quantum affine algebra $U_q(\gh)$.   
However its structure associated with the quantum $Z$-algebras has not yet been discussed so far. 
The purpose of this paper is to address this subject. 
The general theory of the $Z$-algebra 
was studied by Lepowsky and Wilson\cite{LW} and by Gepner\cite{Gepner} in 
the representation theory of affine Lie algebras  and in CFT, respectively.   
Its quantum deformation and application to the representations of $U_q(\gh)$ was partially investigated in \cite{Matsuo, Jing, BoVi, Konno98, JKOPS}.  
A construction of the coset CFT associated with the general 
 $\gh$ was also given \cite{CrRav} in terms of the generalized parafermions. 
We extend these studies to the elliptic algebras $U_{q,p}(\gh)$. In particular,  we 
define a dynamical analogue  $\cZ_k$  of the quantum $Z$-algebras and
 show that the level-$k$ highest weight representations of $U_{q,p}(\gh)$ are 
 realized in terms of 
$\cZ_k$ and the level-$k$ elliptic bosons. It is then shown that the irreducibility of the infinite dimensional 
$U_{q,p}(\gh)$-modules is governed by the $\cZ_k$-modules as in the affine Lie algebra cases\cite{LW}.   

On the other hand,  it was conjectured \cite{Konno98,JKOS} that 
the $U_{q,p}(\gh)$  provides an algebra of the 
screening currents of the deformation of the $W$-algebras associated with the coset 
 $\gh\oplus \gh\supset (\gh)_{{diag}}$ with level $(r-g-1,1)$. For the simply-laced 
 $\gh$, such deformed $W$-algebras have been realized in \cite{FeFr,AKOS,FrRe},   
and in particular for the $\slnh$ case the conjecture has been established 
by an explicit  comparison  of the free field realizations\cite{FeFr,AKOS,KojimaKonno,KK04}. 
However for the non-simply laced $\gh$,  deformation of  the coset type $W$-algebras  has not yet been studied at all. One should note that the coset type 
$W$-algebras associated with the non-simply laced $\gh$ are different from those obtained by the quantum Hamiltonian reduction. See for example \cite{BoSc}. 
We investigate this issue further by giving an explicit realization of the level-1 irreducible highest weight representations of $U_{q,p}(\gh)$ for $\gh=A_l^{(1)}, B_l^{(1)}, D_l^{(1)}, E^{(1)}_6,  E^{(1)}_7, E^{(1)}_8$.  
We show that at least for $A_l^{(1)}$ and $D_l^{(1)}$ the level-1  elliptic currents $e_j(z)$ and $f_j(z)$ coincide with the screening currents of the deformed $W$-algebras obtained in \cite{FeFr,AKOS,FrRe}.  We also show that the irreducible representations  of $U_{q,p}(\gh)$ is naturally decomposed into a direct sum of the irreducible $W$-algebras of the coset type for $\gh=A_l^{(1)}, B_l^{(1)}, D_l^{(1)}$. This suggests in particular an existence of a deformation of Fateev-Lukyanov's $W\! B_l$-algebra \cite{FaLu} as 
the commutant  of the screening operators provided by  the level-1  elliptic currents $e_j(z)$ and $f_j(z)$ of $U_{q,p}(B_l^{(1)})$. 

It is also worth to mention that the coset type $W$-algebras   
 describe a critical behavior of the face type elliptic solvable lattice 
 models\cite{DJKMO, JMO}. 
Correspondingly the $U_{q,p}(\gh)$  provides an algebraic framework to formulate 
 the lattice model itself in the spirit of Jimbo and Miwa\cite{JMBook}. 
This has been established for $\slnh$ in \cite{Konno98,JKOS,KojimaKonno, KKW, Konno08,Konno09} by constructing the $L$-operator and introducing the Hopf 
algebroid structure.  
In order to construct the $L$-operator of $U_{q,p}(\gh)$  and also to get a realization of  a generating function of the deformation of the $W$-algebras,  
 it is crucial to introduce   
new types of elliptic bosons, which we call the fundamental weight type $A^j_m$ and the 
orthonormal basis type $\cE^{\pm j}_m$ distinguishing from the usual ones $\al_{j,m}\ (\al_{j,m}^\vee)$ corresponding to  
 the simple (co-)root and appearing as generators of $U_{q,p}(\gh)$. 
An idea of such bosons {has already appeared} in \cite{FeFr,AKOS,FrRe}. 
 We give an explicit construction of them for 
 $\gh=A_n^{(1)}, B_n^{(1)}, 
 C_n^{(1)}, D_n^{(1)}$. As a check we calculate the commutation relations among  $\cE^{\pm j}_m$ as well as among the elliptic currents $k_{\pm j}(z)$, 
 the generating functions of  $\cE^{\pm j}_m$, and show that they have  
  a universal form. See Theorem \ref{cEcE} and \ref{kk}. 
  
  This paper is organized as follows. In section 2, we define the elliptic algebra $U_{q,p}(\gh)$ as a topological algebra 
 generated by the elliptic Drinfeld generators. This is a new definition of $U_{q,p}(\gh)$ given independently of $U_q(\gh)$ unlike the previous one in Appendix A in \cite{JKOS}. In section 3, we define a quantum dynamical analogue $\cZ_{\hV}$ of Lepowsky and Wilson's $Z$-algebra associated with the level-$k$ $U_{q,p}(\gh)$-module $\hV$ and its universal counterpart $\cZ_k$.  The irreducibility of the level-$k$ highest weight representation of $U_{q,p}(\gh)$ is shown to be governed by the $\cZ_k$-module.  In section 4, we give a simple realization of $\cZ_k$  in terms of the quantum (non-dynamical) $Z$-algebra associated with the level-$k$ $U_q(\gh)$-module and define a standard representation of $U_{q,p}(\gh)$. We provide some level-1 examples of the standard representations and discuss their relation to the deformation of the $W$-algebras.   
  In section 5, we give a construction of the new elliptic bosons of the fundamental weight type and the orthonormal basis type and 
  derive various commutation relations.

\section{Elliptic Algebra $U_{q,p}(\ghbig)$}

\subsection{Definition}
Let   $\gh=X^{(1)}_l$ 
be an untwisted affine {Lie algebra} associated with the generalized Cartan matrix $A=(a_{ij})\ i,j\in \{0\}\cup I,\ I=\{1,\cdots,l\}$. We denote by $B=(b_{ij})$,
 $b_{ij}=d_i a_{ij}$ the symmetrization of $A$. We take $d_i=1\ (i\in I)$ for the simply laced cases, $d_i=1 \ (1\leq i\leq l-1),\ d_l=1/2$ for $B_l^{(1)}$ and 
$d_i=1 \ (1\leq i\leq l-1),\ d_l=2$ for $C_l^{(1)}$. 
Let $q=e^{\hbar}\in \C[[\hbar]]$ and set $q_i=q^{d_i}$. 
Let $p$ be an indeterminate.

Let $\h=\widetilde{\h}\oplus \C d$, $\widetilde{\h}=\bh\oplus\C c$, $\bar{\h}=\oplus_{i\in I}\C h_i$ be the Cartan subalgebra of $\gh$.  
 Define  $\delta, \Lambda_0, \al_i\ (i\in I) \in \h^*$ 
by 
\bea
&&<\al_i,h_j>=a_{j,i}, \ <\delta,d>=1=<\Lambda_0,c>,\lb{pairinghhs}
\ena
the other pairings are 0. We also define $\bar{\Lambda}_i\ (i \in I) \in \h^*$ by 
\be
&&<\bar{\Lambda}_i,h_j>=\delta_{i,j}.
\en
We set $\bh^*=\oplus_{i\in I} \C \bar{\Lambda}_i, $ $\widetilde{\h}^*=\bh^*\oplus \C \Lambda_0$, 
 $\cQ=\oplus_{i\in I}\Z \al_i$ and 
$\cP=\oplus_{i\in I}\Z \bar{\Lambda}_i$. 
 Let $N=l+1$ for $X_l=A_l$, $=l$ for $B_l, C_l, D_l,$  $=7$ for $E_6$, $=8$ for $E_7, E_8$, 
$=3$  for $G_2$, $=4$ for $F_4$ and consider 
the orthonormal basis $\{ \xi_j\ (1\leq j\leq N)\}$  in $\R^N$ 
with the inner product
$( \xi_j, \xi_k )=\delta_{j,k}$.
For $A_l$, we also set   
\bea
&&\bar{\xi}_j=\xi_j-\frac{1}{l+1}\sum_{j=1}^{l+1} \xi_j. 
\ena
We define $\ep_j=\bar{\xi}_j$ for $A_l$ and $=\xi_j$ for other $X_l$.  
The simple roots $\alpha_j$ and the fundamental weights $\bar{\Lambda}_j \ (1\leq j\leq l)$ 
can be expressed as a linear sum of $\ep_j$ 
\cite{Bourbaki, Kac}. We follow Kac's conventions. 
We define $h_{\ep_j}\in \bar{\h} \ (j\in I)$ by 
$<\ep_i,h_{\ep_j}>=(\ep_i,{\ep_j})$ and $h_\alpha \in \bar{\h}$ for $\alpha=\sum_j c_j \ep_j, \ c_j\in \C$ by
 $h_\alpha=\sum_j c_j h_{\ep_j}$.
We regard $\bar{\h}\oplus \bar{\h}^*$ as the Heisenberg 
algebra  by
\bea
&&~[h_{{\epsilon}_j},{\epsilon}_k]
=( {\epsilon}_j,{\epsilon}_k ),\qquad [h_{{\epsilon}_j},h_{{\epsilon}_k}]=0=
[{\epsilon}_j,{\epsilon}_k].\lb{HA1}
\ena
In particular, we have $[h_{j}, \alpha_k]=a_{j k}$.  We also set $h^j=h_{\bar{\Lambda}_j}$. 
 
In order to treat the dynamical shifts in the face type elliptic algebra systematically, 
we  introduce another Heisenberg algebra generated by
$P_{{\alpha}}$ and $ Q_{{\beta}}\ 
({\alpha}, {\beta} \in \bar{\h}^*)$ satisfying the 
commutation relations
\begin{eqnarray}
&&[P_{{\epsilon}_j}, Q_{{\epsilon}_k}]=
( {\epsilon}_j, {\epsilon}_k ), \qquad 
[P_{{\epsilon}_j}, P_{{\epsilon}_k}]=0=
[Q_{{\epsilon}_j}, Q_{{\epsilon}_k}].\lb{HA2}
\end{eqnarray}
{We also set}
\begin{eqnarray}
&&[P_{{\epsilon}_j}, \al]=
[Q_{{\epsilon}_j}, \al]=0,\qquad  
[P_{{\epsilon}_j}, U(\gh)]=
[Q_{{\epsilon}_j}, U(\gh)]=0\lb{HA4} 
\end{eqnarray}
where
$P_\alpha=\sum_j c_j P_{\bep_j}$ for $\alpha=\sum_j c_j \bep_j$.
We set $P_{\bh}=\oplus_{j\in I}\C P_{\bep_j}, Q_{\bh}=\oplus_{j\in I}\C Q_{\bep_j}$  $P_{j}=P_{\alpha_j^\vee
}, P^j=P_{\bar{\Lambda}_j}$ and 
$Q_{j}=Q_{\alpha_j}, Q^j=Q_{\bar{\Lambda}_j^\vee}$. Here $\alpha _{j}^\vee=2\alpha _{j}/(\al_j,\al_j)$.

For the abelian group 
$\cR_Q= \sum_{j=1}^N\Z Q_{\al_j}$, 
we denote by $\C[\cR_Q]$ the group algebra over $\C$ of $\cR_Q$. 
We denote by $e^{\al}$ the element of $\C[\cR_Q]$ corresponding to $\al\in \cR_Q$. 
These $e^\al$ satisfy $e^\al e^\beta=e^{\al+\beta}$
 and $(e^\al)^{-1}=e^{-\al}$. 
In particular, $e^0=1$ is the identity element. 

Now let us set $H=\widetilde{\h}\oplus P_{\bh}=\sum_{j}\C(P_{\bep_j}+h_{\bep_j})+\sum_j\C P_{\bep_j}+\C c$ and denote its dual space by $H^*=
\widetilde{\h}^*\oplus Q_{\bh}$. We define the paring by \eqref{pairinghhs},  
$<Q_\al,P_\beta>=(\al,\beta)$ and  $<Q_\al,h_\beta>=<Q_\al,c>=<Q_\al,d>=0=<\al,P_\beta>=<\delta,P_\beta>=<\Lambda_0,P_\beta>$ . We define $\K=\cM_{H^*}$ to be the field of meromorphic functions on $H^*$.  We regard a function of $P+h=\sum_ja_j(P_{\bep_j}+h_{\bep_j})$, $P=\sum_j b_jP_{\bep_j}$ and $c$,  $\hf=f(P+h,P, c)$, as an element in $\K$ by   
 $\hf(\mu)=f(<\mu,P+h>, <\mu,P>, <\mu,c>)$ for  $\mu\in H^*$. 

We use the following notations.
\be 
&&[n]=\frac{q^n-q^{-n}}{q-q^{-1}}, \quad [n]_i=\frac{q_i^n-q_i^{-n}}{q_i-q_i^{-1}}, \qquad [n)_j=\frac{q^n-q^{-n}}{q_j-q_j^{-1}},\\
&&[n]_i!=[n]_i[n-1]_i\cdots [1]_i,\quad \left[\mmatrix{m\cr n\cr}\right]_i=\frac{[m]_i!}{[n]_i![m-n]_i!},\\
&&(x;q)_\infty=\prod_{n=0}^\infty(1-x q^n),\quad (x;q,t)_\infty=\prod_{n,m=0}^\infty(1-x q^n t^m),\quad
\Theta_p(z)=(z;p)_{\infty}(p/z;p)_\infty(p;p)_\infty. 
\en

\begin{dfn}\lb{defUqp}
{The elliptic algebra} $U_{q,p}(\gh)$ is a topological algebra over $\K[[p]]$ generated by $\cM_{H^*}$,  $e_{j,m}, f_{j,m}, \al^\vee_{j,n}, K^\pm_j$,  
$(j\in I, m\in \Z, n\in \Z_{\not=0})$, ${\hd}$ and the central element $c$. 
We assume $K^\pm_{j}$ are invertible and set 
\be
&&e_j(z)=\sum_{m\in \Z} e_{j,m}z^{-m}
,\quad f_j(z)=\sum_{m\in \Z} f_{j,m}z^{-m},\lb{deffn}\\
&&{\psi}_j^+(
q^{-\frac{c}{2}}
z)=K^+_{j}\exp\left(-(q_j-q_j^{-1})\sum_{n>0}\frac{\al^\vee_{j,-n}}{1-p^n}z^n\right)
\exp\left((q_j-q_j^{-1})\sum_{n>0}\frac{p^n\al^\vee_{j,n}}{1-p^n}z^{-n}\right),\\
&&{\psi}_j^-(q^{\frac{c}{2}}z)=K^-_{j} \exp\left(-(q_j-q_j^{-1})\sum_{n>0}\frac{p^n\al^\vee_{j,-n}}{1-p^n}z^n\right)
\exp\left((q_j-q_j^{-1})\sum_{n>0}\frac{\al^\vee_{j,n}}{1-p^n}z^{-n}\right). 
\en
Note that  $\psi^\pm_j(z)$ are formal Laurent series in $z$, whose coefficients 
are well defined in the $p$-adic topology. 
We call $e_j(z), f_j(z), \psi^\pm_j(z)$ the elliptic currents. 
The  defining relations are as follows.  For $g(P), g(P+h)\in \cM_{H^*}$, 
\bea
&&g({P+h})e_j(z)=e_j(z)g({P+h}),\quad g({P})e_j(z)=e_j(z)g(P-<Q_{\al_j},P>),\lb{ge}\\
&&g({P+h})f_j(z)=f_j(z)g(P+h-<{\al_j},P+h>),\quad g({P})f_j(z)=f_j(z)g(P),\lb{gf}\\
&&[g(P), \al^\vee_{i,m}]=[g(P+h),\al^\vee_{i,n}]=0,\lb{gboson}\\ 
&&g({P})K^\pm_j=K^\pm_jg(P-<Q_{\al_j},P>),\\
&&\ g({P+h})K^\pm_j=K^\pm_jg(P+h-<Q_{\al_j},P>),\lb{gKpm}
\\
&&[\hd, g(P)]=[\hd, g(P+h)]=0,
\quad\lb{dg}\\
&& [\hd, \al^\vee_{j,n}]=n\al^\vee_{j,n},\quad [\hd, e_j(z)]=-z\frac{\partial}{\partial z}e_j(z), \quad  [\hd, f_j(z)]=-z\frac{\partial}{\partial z}f_j(z),\quad \lb{dedf}\\
&&K_i^{\pm}e_j(z)=q_i^{\mp a_{ij}}e_j(z)K_i^{\pm},\quad 
K_i^{\pm}f_j(z)=q_i^{\pm a_{ij}}f_j(z)K_i^{\pm},
\\
&&[\al^\vee_{i,m},\al^\vee_{j,n}]=\delta_{m+n,0}\frac{[a_{ij}m]_i
[cm)_j
}{m}
\frac{1-p^m}{1-p^{*m}}
q^{-cm}
,\lb{ellboson}\\
&&
[\al^\vee_{i,m},e_j(z)]=\frac{[a_{ij}m]_i}{m}\frac{1-p^m}{1-p^{*m}}
q^{-cm}z^m e_j(z),
\lb{bosonve}\\
&&
[\al^\vee_{i,m},f_j(z)]=-\frac{[a_{ij}m]_i}{m}z^m f_j(z)
,\lb{bosonvf}\\
&&
z_1 \frac{(q^{b_{ij}}z_2/z_1;p^*)_\infty}{(p^*q^{-b_{ij}}z_2/z_1;p^*)_\infty}e_i(z_1)e_j(z_2)=
-z_2 \frac{(q^{b_{ij}}z_1/z_2;p^*)_\infty}{(p^*q^{-b_{ij}}z_1/z_2;p^*)_\infty}e_j(z_2)e_i(z_1),\lb{ee}\\
&&
z_1 \frac{(q^{-b_{ij}}z_2/z_1;p)_\infty}{(pq^{b_{ij}}z_2/z_1;p)_\infty}f_i(z_1)f_j(z_2)=
-z_2 \frac{(q^{-b_{ij}}z_1/z_2;p)_\infty}{(pq^{b_{ij}}z_1/z_2;p)_\infty}f_j(z_2)f_i(z_1),\lb{ff}\\
&&[e_i(z_1),f_j(z_2)]=\frac{\delta_{i,j}}{q_i-q_i^{-1}}
\left(\delta(
q^{-c}
z_1/z_2)
\psi_j^-(
q^{\frac{c}{2}}
z_2)-
\delta(
q^c
z_1/z_2)
\psi_j^+(
q^{-\frac{c}{2}}
z_2)
\right),\lb{eifj}
\ena
\bea
&&\sum_{\sigma\in S_{a}}\prod_{1\leq m< k \leq a}
\frac{(p^*q^{2}{z_{\sigma(k)}}/{z_{\sigma(m)}}; p^*)_{\infty}}
{(p^*q^{-2}{z_{\sigma(k)}}/{z_{\sigma(m)}}; p^*)_{\infty}}\nn\\
&&\quad\times\sum_{s=0}^{a}(-1)^{s}
\left[\begin{array}{c}
a\cr
s\cr
\end{array}\right]_{i}\prod_{1\leq m \leq s}\frac{(p^*q^{b_{ij}}{w}/{z_{\sigma(m)}}; p^*)_{\infty}}
{(p^*q^{-b_{ij}}{w}/{z_{\sigma(m)}}; p^*)_{\infty}}
\prod_{s+1\leq m \leq a}\frac{(p^*q^{b_{ij}}{z_{\sigma(m)}}/{w}; p^*)_{\infty}}
{(p^*q^{-b_{ij}}{z_{\sigma(m)}}/{w}; p^*)_{\infty}}
\nn\\
&&\quad\times
e_{i}(z_{\sigma(1)})\cdots e_{i}(z_{\sigma(s)})e_{j}(w)e_{i}(z_{\sigma(s+1)}) \cdots 
e_{i}(z_{\sigma(a)})=0,\label{serree}
\ena
\bea
&&\sum_{\sigma\in S_{a}}\prod_{1\leq m< k \leq a
}
\frac{(pq^{-2}{z_{\sigma(k)}}/{z_{\sigma(m)}}; p)_{\infty}}
{(pq^{2}{z_{\sigma(k)}}/{z_{\sigma(m)}}; p)_{\infty}}\nn\\
&&\quad\times\sum_{s=0}^{a}(-1)^{s}\left[\begin{array}{c}
a\cr
s\cr
\end{array}\right]_{i}\prod_{1\leq  m \leq s}
\frac{(pq^{-b_{ij}}{w}/{z_{\sigma(m)}}; p)_{\infty}}{(pq^{b_{ij}}{w}/{z_{\sigma(m)}}; p)_{\infty}}
\prod_{s+1\leq m \leq a}\frac{(pq^{-b_{ij}}{z_{\sigma(m)}}/{w}; p)_{\infty}}{(pq^{b_{ij}}{z_{\sigma(m)}}/{w}; p)_{\infty}}\nn\\
&&\times f_{i}(z_{\sigma(1)})\cdots f_{i}(z_{\sigma(s)})f_{j}(w)f_{i}(z_{\sigma(s+1)}) \cdots f_{i}(z_{\sigma(a)})=0
\quad(i\neq j, a=1-a_{ij}),\label{serref}
\end{eqnarray}
where $p^*=p
q^{-2c}
$ and $\delta(z)=\sum_{n\in \Z}z^n$.  
We also denote by ${U}'_{q,p}(\gh)$ the subalgebra obtained by removing ${d}$. 

\end{dfn}
We treat the relations \eqref{dedf}, \eqref{bosonve}-\eqref{serref} as formal Laurent series in $z, w$ and $z_j$'s. 
In  each term of \eqref{ee}-\eqref{serref}, the expansion direction of the structure function given by a ratio of  infinite products is 
chosen according to the order of the accompanied product of the
 elliptic currents. For example, in the l.h.s of \eqref{ee}, $\frac{(q^{b_{ij}}z_2/z_1;p^*)_\infty}{(p^*q^{-b_{ij}}z_2/z_1;p^*)_\infty}$ should be  expanded in $z_2/z_1$, 
 whereas in the r.h.s $\frac{(q^{b_{ij}}z_1/z_2;p^*)_\infty}{(p^*q^{-b_{ij}}z_1/z_2;p^*)_\infty}$ should be expanded in $z_1/z_2$.   
 In each term in \eqref{serree}, the coefficient function is expanded in  $z_{\sigma{(k)}}/z_{\sigma{(m)}}\ (m<k)$, 
 $w/z_{\sigma(m)}\ (m\leq s)$ and $z_{\sigma(m)}/w\ (m\geq s+1)$. 
All the coefficients in $z_j$'s are well defined in the $p$-adic topology.  

\noindent
{\it Remark.}\ 
In \cite{Konno98, JKOS, KKW, Konno08}, assuming that $q$ is a transcendental complex number  satisfying $|q|<1$, we wrote \eqref{ee}, \eqref{ff} as
\be
&&z_1 \Theta_{p^*}(q^{b_{ij}}z_2/z_1)e_i(z_1)e_j(z_2)=
-z_2 \Theta_{p^*}(q^{b_{ij}}z_1/z_2)e_j(z_2)e_i(z_1),\lb{eetheta}\\
&&z_1 \Theta_{p}(q^{-b_{ij}}z_2/z_1)f_i(z_1)f_j(z_2)=
-z_2 \Theta_{p}(q^{-b_{ij}}z_1/z_2)f_j(z_2)f_i(z_1),\lb{fftheta}
\en
in the sense of analytic continuation.

Let $U_q(\gh)$ be the quantum affine algebra associated with $\gh$ in the 
Drinfeld realization\cite{Dr}. 
See Appendix \ref{Uqgh}.  $U_{q,p}(\gh)$ is a natural face type ( i.e. dynamical) 
elliptic deformation of $U_q(\gh)$ in the following sense. 
\begin{thm}
\be
&&{U_{q,p}(\gh)/pU_{q,p}(\gh)\cong (\K\otimes_{\C}U_q(\gh))\sharp \C[\cR_Q]}. 
\en
Here the smash product $\sharp$ is defined as follows.
\be
&&g(P,P+h)x\otimes e^\al \cdot f(P,P+h)y\otimes e^\beta\\
&&\quad= g(P,P+h)f(P-<\al,P>,P+h-<\al+\wt x,P+h>)xy\otimes e^{\al+\beta}
\en
where $\wt{x}\in \bar{\h}^*$ s.t. $q^h x q^{-h}=q^{<\wt{x},h>}x$ for $x, y\in U_q(\gh), f(P), g(P)\in \K, e^\al, e^\beta\in \C[\cR_Q]$.  
\end{thm}
\noindent
{\it Proof.}\ At $p=0$, the relations for $\al^\vee_{j,m}, e_j(z), f_j(z)$\ \eqref{dedf}-\eqref{serref} coincide with those for $a^\vee_{i,m}, x^+_j(z), x^-_j(z)$\    \eqref{aaUq}-\eqref{serrefUq} 
of $U_q(\gh)$. Therefore from \eqref{ge}-\eqref{gKpm}, one has the isomorphism
\be
&&e_j(z)\mapsto x^+_j(z)e^{-Q_{\al_j}},\quad f_j(z)\mapsto x^-_j(z), \quad 
K^\pm_j\mapsto q_j^{\mp h_j}e^{-Q_{\al_j}},\quad \al^\vee_{j,m}\mapsto a^\vee_{j,m}\quad\mbox{mod}\ pU_{q,p}(\gh). 
\en 

\subsection{$H$-algebra $U_{q,p}(\ghbig)$}
Let $\cA$ be a complex associative algebra, $\cH$ be a finite dimensional commutative subalgebra of 
$\cA$, and $\cM_{\cH^*}$ be the 
field of meromorphic functions on $\cH^*$ the dual space of $\cH$. 

\begin{dfn}[$\cH$-algebra] 
An $\cH$-algebra is an associative algebra $\cA$ with 1, which is bigraded over 
$\cH^*$, $\ds{\cA=\bigoplus_{\alpha,\beta\in \cH^*} \cA_{\al\beta}}$, and equipped with two 
algebra embeddings $\mu_l, \mu_r : \cM_{\cH^*}\to \cA_{00}$ (the left and right moment maps), such that 
\be
\mu_l(\hf)a=a \mu_l(T_\al \hf), \quad \mu_r(\hf)a=a \mu_r(T_\beta \hf), \qquad 
a\in \cA_{\al\beta},\ \hf\in \cM_{\cH^*},
\en
where $T_\al$ denotes the automorphism $(T_\al \hf)(\la)=\hf(\la+\al)$ of $\cM_{\cH^*}$.
\end{dfn}

\begin{prop}
$U=U_{q,p}(\gh)$  is {an $H$-algebra} by 
\be
&&U=\bigoplus_{\al,\beta\in H^*}U_{\al,\beta}\\
&&U_{\al,\beta}=\left\{x\in U \left|\ q^{P+h}x q^{-(P+h)}=q^{<\al,P+h>}x,\quad q^{P}x q^{-P}=q^{<\beta,P>}x,\ \forall P+h, P\in H\right.\right\}
\en
and $\mu_l, \mu_r : \K \to U_{0,0}$ defined by 
\be
&&\mu_l(\hf)=f(P+h,p)\in \K[[p]],\qquad \mu_r(\hf)=f(P,p^*)\in \K[[p]].
\en
\end{prop}

\subsection{Dynamical Representations}
Let us consider a vector space $\hV$ over $\K$, which is  
${H}$-diagonalizable, i.e.  
\be
&&\hV=\bigoplus_{\la,\mu\in {H}^*}\hV_{\la,\mu},\ \hV_{\la,\mu}=\{ v\in \hV\ |\ q^{P+h}\cdot v=q^{<\la,P+h>} v,\ q^{P}\cdot v=q^{<\mu,P>} v\ \forall 
P+h, P\in 
{H}\}.
\en
Let us define the $H$-algebra $\cD_{H,\hV}$ of the $\C$-linear operators on $\hV$ by
\be
&&\cD_{H,\hV}=\bigoplus_{\al,\beta\in {H}^*}(\cD_{H,\hV})_{\al\beta},\\
&&\hspace*{-10mm}(\cD_{H,\hV})_{\al\beta}=
\left\{\ X\in \End_{\C}\hV\ \left|\ 
\mmatrix{ f(P+h)X=X f(P+h+<\alpha,P+h>){},\cr 
f(P)X=X f(P+<\beta,P>){,}\cr
 f(P), f(P+h)\in \K,\ X\cdot\hV_{\la,\mu}\subseteq 
 \hV_{\la+\al,\mu+\beta}\cr}  
 \right.\right\},\\
&&\mu_l^{\cD_{H,\hV}}(\widehat{f})v=f(<\la,P+h>,p)v,\quad 
\mu_r^{\cD_{H,\hV}}(\widehat{f})v=f(<\mu,P>,p^*)v,\quad \widehat{f}\in {\cM}_{H^*},
\ v\in \hV_{\la,\mu}.
\en
\begin{dfn}\lb{DRep}
We define a dynamical representation of $U_{q,p}(\gh)$ on $\hV$ to be  
 an $H$-algebra homomorphism ${\pi}: U_{q,p}(\gh) 
 \to \cD_{H,\hV}$. By the action $\pi$ of $U_{q,p}(\gh)$ we regard $\hV$ as a 
$U_{q,p}(\gh)$-module. 
\end{dfn}
%
\begin{dfn}
For $k\in \C$, we say that a $U_{q,p}(\gh)$-module has  level $k$ if $c$ act 
as the scalar $k$ on it.  
\end{dfn}
\noindent
{\it Remark.} For the level-0 representations,  Definition \ref{DRep} is essentially the same as in \cite{EV}, by identifying $P$ and $P+h$ with $\la$ and $\la-\gamma h$, 
respectively. This definition is valid also for the non-zero level cases\cite{Konno09}. 
\begin{dfn}
For $\omega\in \C$, we set 
\be
&&\hV_\omega=\{v\in \hV\ |\ \hd\cdot v=\omega v\ \}
\en
and we call $\hV_\omega$ the space of elements homogeneous of degree $\omega$. 
We also say that $X\in \cD_{H,\hV}$ is homogeneous of degree $\omega\in \C$ if 
\be
&&[\hd, X]=\omega X
\en
and denote by $(\cD_{H,\hV})_\omega$ the space of all {endomorphisms} homogeneous of degree $\omega$. 
\end{dfn}

\begin{dfn}
Let ${\cal H}$, ${\cal N}_+, {\cal N}_-$ be the subalgebras of 
$U_{q,p}(\gh)$ 
generated by 
$c, d, 
K^\pm_{i}\ (i\in I)$, by $\al^\vee_{i,n}\ (i\in I, n\in \Z_{>0})$,  $e_{i, n}\ (i\in I, n\in \Z_{\geq 0})$  
$f_{i, n}\ (i\in I, n\in \Z_{>0})$ and by $\al^\vee_{i,-n}\ (i\in I, n\in \Z_{>0}),\ e_{i, -n}\ (i\in I, n\in \Z_{> 0}),\ f_{i, -n}\ (i\in I, n\in \Z_{\geq 0})$, respectively.   
\end{dfn}

\begin{dfn}
For $k\in\C$, $\la\in \h^*$ and $\mu\in H^*$, 
a (dynamical) $U_{q,p}(\gh)$-module $\hV(\la,\mu)$ is called the 
level-$k$ highest weight module with the highest weight $(\la,\mu)$, if there exists a vector 
$v\in \hV(\la,\mu)$ such that
\be
&&\hV(\la,\mu)=U_{q,p}(\gh)\cdot v,\qquad \cN_+\cdot v=0,\\
&&c\cdot v=kv, 
\quad  f({P})\cdot v =f({<\mu,P>})v,\quad f({P+h})\cdot v =f({<\la,P+h>})v.
\en
\end{dfn}

We define the category $\cC_k$ in the analogous way to the classical affine Lie algebra case\cite{LW}. 
\begin{dfn}
For $k\in \C$, $\cC_k$ is the full subcategory of the category of $U_{q,p}(\gh)$-modules 
consisting of those modules $\hV$ such that 
\begin{itemize}
\item[(i)] $\hV$ has level $k$
\item[(ii)] $\ds{\hV=\bigsqcup_{\omega\in \C} \hV_\omega}$
\item[(iii)] For every $\omega \in \C$, there exists $n_0\in \N$ such that for all 
$n>n_0$, $\hV_{\omega+n}=0$.   
\end{itemize}
\end{dfn}

Since $\hpi \cN_+\subset \bigsqcup_{n\in \Z_{\geq 0}}(\cD_{H,\hV})_n$, any level-$k$ highest weight $U_{q,p}(\gh)$-modules belong to $\cC_k$.

\section{The Dynamical Quantum $Z$-Algebras}

In this section we introduce a quantum and  dynamical analogue $\cZ_k$ of Lepowsky-Wilson's $Z$-algebra associated with the level-$k$ 
$U_{q,p}(\gh)$-modules
and define a category $\fD_k$ of the $\cZ_k$-modules. 
Each representation of $\cZ_k$ in $\fD_k$ turns out to be a dynamical analogue of the 
quantum $Z$-algebra derived by Jing \cite{Jing} from the level-$k$ representation in 
the  $U_{q}(\gh)$ counterpart of $\fD_k$. See sec.4.1.
We also provide the Serre relations (\ref{z5}) which are not written in \cite{Jing} explicitly. 

\subsection{The Heisenberg algebra $U_{q,p}(\cH)$}
Let $U_{q,p}(\cH)$ be the subalgebra of $U_{q,p}(\gh)$ 
 generated by  $\al^\vee_{i,n}\ (i\in I, n\in \Z_{\not=0})$ and 
$c$.   
It is convenient to introduce the simple root type generators $\al_{j,m}$ and $\al'_{j,m}$ defined by $\al_{j,m}=[d_j]\al^\vee_{j,m}$ 
and $\ds{\alpha'_{j,m}=\frac{1-p^{*m}}{1-p^m}q^{cm} \alpha_{j,m},\ (j\in I, m\neq 0)}$. 
From \eqref{ellboson}, \eqref{bosonve}, \eqref{bosonvf},  we have 
\bea
&&[\al_{i,m},\al_{j,n}]=\frac{[b_{ij}m][cm]}{m}\frac{1-p^m}{1-p^{*m}}q^{-cm}\delta_{m+n,0},\lb{alal}\\
&&[\al'_{i,m},\al'_{j,n}]=\frac{[b_{ij}m][cm]}{m}\frac{1-p^{*m}}{1-p^{m}}q^{cm}\delta_{m+n,0},\\
&&[\al_{i,m},\al'_{j,n}]=\frac{[b_{ij}m][cm]}{m}\delta_{m+n,0},
\ena
\bea
&&[\al_{i,m},e_j(z)]=\frac{[b_{ij}m]}{m}\frac{1-p^m}{1-p^{*m}}q^{-cm}z^m e_j(z),
\lb{bosone}\\
&&[\al'_{i,m},f_j(z)]=-\frac{[b_{ij}m]}{m}\frac{1-p^{*m}}{1-p^{m}}q^{cm}z^m f_j(z).\lb{bosonf}
\ena

Let 
$U_{q,p}(\cH_+)$ (resp. $U_{q,p}({\cH}_-)$) be the commutative
 subalgebras of  $U_{q,p}({{\cH}})$ generated by $\{c, \alpha_{i,n}(i\in I, n\in\Z_{> 0})\}$ 
(resp. $\{\alpha_{i,-n}(i\in I, n\in\Z_{> 0})\}$). 
We have
\be
U_{q,p}(\cH)=U_{q,p}(\cH_-)U_{q,p}(\cH_+).
\en

Let  $\C 1_{k}$ be the one-dimensional 
 $U_{q,p}(\cH_+)$-module generated by the vacuum vector $1_k$ defined by  
\be
&&c\cdot 1_k=k1_k
\qquad \alpha_{i,n}\cdot 1_{k}=0 \qquad(n> 0). 
\en
Then we have the induced $U_{q,p}({\cH})$-module
\be
\F_{\al,k}=U_{q,p}(\cH)\otimes_{U_{q,p}({\cH}_+)}\C1_{k}.
\en
We identify $\F_{\al,k}$ with a polynomial ring $\C[\al_{i,-m}\ (i\in I,m>0)]$ 
by  
\be
&&c\cdot u= ku,\quad  
\alpha_{i,-n}\cdot u= \alpha_{i,-n}u,\qquad \\
&&\alpha_{i,n}\cdot u=\sum_{j} \frac{[b_{ij}n][kn]}{n}\frac{1-p^n}{1-p^{*n}}q^{-kn}\frac{\partial}{\partial\alpha_{j,-n}}u\qquad(n>0)
\en
for $u\in \C[\al_{i,-m}\ (i\in I,m>0)]$. 

\subsection{The dynamical quantum $Z$-algebra $\cZ_{\hV}$}
Let $k\in\C^\times$ and  $(\hV,\hpi)
\in \sC_k$.  We call $\hpi U_{q,p}(\cH)\subset (D_{H,\hV})_{00}$ the level-$k$  Heisenberg algebra.  
We define the following vertex operators in $(D_{H,\hV})_{00}[[z,z^{-1}]]$.
\be
&&E^\pm(\al_j,z)=\exp\left(\pm\sum_{n> 0}\frac{\hpi(\alpha_{j,\pm n})}{[kn]}z^{\mp n}\right),
\qquad E^\pm(\al'_j,z)=\exp\left(\mp\sum_{n> 0}\frac{\hpi(\alpha'_{j,\pm n})}{[kn]}z^{\mp n}\right).
\en
These satisfy the following relations.
\begin{prop} 
\bea
&&E^+(\al_i,z)E^-(\al_j,w)=\frac{(q^{-b_{ij}+2k}w/z;q^{2k})_\infty(q^{-b_{ij}}w/z;p^*)_\infty}{(q^{b_{ij}+2k}w/z;q^{2k})_\infty
(q^{b_{ij}}w/z;p^*)_\infty}E^-(\al_j,w)E^+(\al_i,z),\lb{EpEp}\\
&&E^+(\al'_i,z)E^-(\al'_j,w)=\frac{(q^{-b_{ij}}w/z;q^{2k})_\infty(q^{b_{ij}}w/z;p)_\infty}{(q^{b_{ij}}w/z;q^{2k})_\infty
(q^{-b_{ij}}w/z;p)_\infty}E^-(\al'_j,w)E^+(\al'_i,z),\\
&&E^+(\al_i,z)E^-(\al'_j,w)=\frac{(q^{b_{ij}+k}w/z;q^{2k})_\infty}{(q^{-b_{ij}+k}w/z;q^{2k})_\infty}E^-(\al'_j,w)E^+(\al_i,z),\\
&&E^+(\al'_i,z)E^-(\al_j,w)=\frac{(q^{b_{ij}+k}w/z;q^{2k})_\infty}{(q^{-b_{ij}+k}w/z;q^{2k})_\infty}E^-(\al_j,w)E^+(\al'_i,z),
\ena
\bea
&&E^\pm(\al_i,z)e_j(w)=\frac{(q^{\pm b_{ij}+2k}(w/z)^{\pm 1};q^{2k})_\infty(q^{\pm b_{ij}}(w/z)^{\pm 1};p^*)_\infty}{(q^{\mp b_{ij}+2k}(w/z)^{\pm 1};q^{2k})_\infty
(q^{\mp b_{ij}}(w/z)^{\pm 1};p^*)_\infty}e_j(w)E^\pm(\al_i,z),\\
&&E^\pm(\al'_i,z)f_j(w)=\frac{(q^{\pm b_{ij}}(w/z)^{\pm 1};q^{2k})_\infty(q^{\pm b_{ij}}(w/z)^{\pm 1};p)_\infty}{(q^{\mp b_{ij}}(w/z)^{\pm 1};q^{2k})_\infty
(q^{\mp b_{ij}}(w/z)^{\pm 1};p)_\infty}f_j(w)E^\pm(\al'_i,z),\\
&&E^\pm(\al'_i,z)e_j(w)=\frac{(q^{\mp b_{ij}+k}(w/z)^{\pm 1};q^{2k})_\infty
}{(q^{\pm b_{ij}+k}(w/z)^{\pm 1};q^{2k})_\infty
}e_j(w)E^\pm(\al'_i,z),\\
&&E^\pm(\al_i,z)f_j(w)=\frac{(q^{\mp b_{ij}+k}(w/z)^{\pm 1};q^{2k})_\infty
}{(q^{\pm b_{ij}+k}(w/z)^{\pm 1};q^{2k})_\infty
}f_j(w)E^\pm(\al_i,z).
\ena
\end{prop}

\begin{dfn}
We define $\cZ^\pm_j(z;\hV)\in \cD_{H,\hV}[[z,z^{-1}]]$ by
\bea
&&\cZ^+_j(z;\hV):=E^-(\alpha_{j},z)\hpi(e_{j}(z))E^+(\al_j,z)
,\lb{zpe}\\
&&\cZ^-_j(z;\hV):=E^-(\alpha'_{j},z)\hpi(f_{j}(z))E^+(\al'_j,z).\lb{zmf}
\ena
for $j\in I$ and call them  the dynamical quantum $Z$ operators associated with $(\hV,\pi)\in \fC_k$. 
\end{dfn}
Note that due to the truncation property of the grading of $\hV\in \sC_k$ w.r.t $\hd$,   
$\cZ^\pm_j(z;\hV)$ are well defined i.e. the coefficients $\cZ_{j,n}^\pm(\hV)$ of  
$\cZ^\pm_j(z;\hV)=\sum_{n\in \Z}\cZ^\pm_{j,n}(\hV) z^{-n}$ in $z$ are well defined elements in   
 $(\cD_{H,\hV})_n$ for all $n\in \Z$. For the sake of simplicity of the presentation, 
we often drop $\hpi$ to denote the elements in $\cD_{H,\hV}$.    

From the defining relations of $U_{q,p}(\gh)$, we obtain the following relations of the dynamical quantum $Z$ operators. 
\begin{thm}\lb{zz}
\bea
&&\hspace{-1.6cm}g({P+h})\cZ^+_i(z;\hV)=\cZ^+_i(z;\hV)g{(P+h)},\quad 
g({P})\cZ^+_i(z;\hV)=\cZ^+_i(z;\hV)g(P-<Q_{\al_i},P>),
\lb{hZPZp}
\\
&&\hspace{-1.6cm}g({P+h})\cZ^-_i(z;\hV)=\cZ^-_i(z;\hV)g(P+h-<{\al_i},P+h>),\quad 
g({P})\cZ^-_i(z;\hV)=\cZ^-_i(z;\hV)g(P),\lb{hZPZm}\\
&&[\hd, \cZ^\pm_j(z;\hV)]=-z\frac{\partial}{\partial z}\cZ^\pm_j(z;\hV),\lb{dZ}\\
&&[\alpha_{i,m},\cZ^\pm_j(w;\hV)]=0,\lb{bosonZ}\\
&&K_i^\pm \cZ^+_{j}(z;\hV)=q^{\mp b_{ij}}\cZ^+_{j}(z;\hV)K_i^\pm,\quad 
K_i^\pm \cZ^-_{j}(z;\hV)=q^{\pm b_{ij}}\cZ^-_{j}(z;\hV)K_i^\pm,
\ena
\bea
&&\hspace{-0.5cm}\label{z1}z\frac{(q^{-b_{ij}}{w}/{z};q^{2k})_{\infty}}{(q^{b_{ij}+2k}{w}/{z};q^{2k})_{\infty}}
\cZ^\pm_i(z;\hV)\cZ^\pm_j(w;\hV)=
-w\frac{(q^{-b_{ij}}{z}/{w};q^{2k})_{\infty}}{(q^{b_{ij}+2k}{z}/{w};q^{2k})_{\infty}}
\cZ^\pm_j(w;\hV)\cZ^\pm_i(z;\hV),\\
&&\frac{(q^{b_{ij}+k}{w}/{z};q^{2k})_{\infty}}{(q^{-b_{ij}+k}{w}/{z};q^{2k})_{\infty}}
\cZ^+_i(z;\hV)\cZ^-_j(w;\hV)
-\frac{(q^{b_{ij}+k}{z}/{w};q^{2k})_{\infty}}{(q^{-b_{ij}+k}{z}/{w};q^{2k})_{\infty}}
\cZ^-_j(w;\hV)\cZ^+_i(z;\hV)\nn \\
&&\label{z2}\qquad\qquad\qquad\qquad=\frac{\delta_{ij}}{q_{i}-q_{i}^{-1}}
\left(K_i^-\delta\bigl(q^{-k}{z}/{w}\bigr)
-K^+_i\delta\bigl(q^{k}{z}/{w}\bigr)
\right),
\ena
\bea
&&\sum_{\sigma\in S_{a}}\prod_{1\leq m< l \leq a} 
\frac{(q^{2+k\mp k}{z_{\sigma(l)}}/{z_{\sigma(m)}}; q^{2k})_{\infty}}
{(q^{-2+k\mp k}{z_{\sigma(l)}}/{z_{\sigma(m)}}; q^{2k})_{\infty}}
\nn\\
&&\quad\times\sum_{s=0}^{a}(-1)^{s}\left[\begin{array}{c}
a\cr
s\cr
\end{array}\right]_{i}\prod_{1\leq m\leq s}
\frac{(q^{-b_{ij}+k\mp k}{w}/{z_{\sigma(m)}};q^{2k})_{\infty}}
{(q^{b_{ij}+k\mp k}{w}/{z_{\sigma(m)}}; q^{2k})_{\infty}}
\prod_{s+1\leq m\leq a}\frac{(q^{-b_{ij}+k\mp k}{z_{\sigma(m)}}/{w}; q^{2k})_{\infty}}
{(q^{b_{ij}+k\mp k}{z_{\sigma(m)}}/{w};q^{2k})_{\infty}}
\nn\\
&&\label{z5}\quad \times \cZ^{\pm}_{i}(z_{\sigma(1)};\hV)\cdots \cZ^{\pm}_{i}(z_{\sigma(s)};\hV)\cZ^{\pm}_{j}(w;\hV)\cZ^{\pm}_{i}
(z_{\sigma(s+1)};\hV)\cdots \cZ^{\pm}_{i}(z_{\sigma(a)};\hV)=0\nn\\
&&\qquad\qquad(i\not= j,\ a=1-a_{ij}). \lb{SerreZ}
\ena 
\end{thm}

\noindent
{\it Proof.}\ 
The relations \eqref{hZPZm} and \eqref{dZ} follow from \eqref{ge}-\eqref{gKpm} and \eqref{dedf}, respectively.  
Let us show the relation (\ref{bosonZ}). 
For $m>0$, we have 
\be
[\al_{i,m},\cZ^+_j(z;\hV)]&=&[\al_{i,m},E^-(\al_j,z)]e_{j}(z)E^+(\al_j,z)
+E^-(\al_j,z)[\al_{i,m},e_j(z)]E^+(\al_j,z).
\en
This vanishes due to \eqref{bosone} and 
\be
&&[\al_{i,m},E^-(\al_j,z)]=-\frac{[b_{ij}m]}{m}\frac{1-p^m}{1-p^{*m}}q^{-km}z^{m},
\en
where $p^*=pq^{-2k}$. 
Similarly, $[\al_{i,m},\cZ^-_j(z;\hV)]=0$ follows from \eqref{bosonf} and 
\be
&&[\al'_{i,m},E^-(\al'_j,z)]=\frac{[b_{ij}m]}{m}\frac{1-p^{*m}}{1-p^m}q^{km}z^{m}.
\en
The case $m<0$ can be proved in a similar way.  

The relation \eqref{z1} follows from
\be
&&\cZ^+_i(z;\hV)\cZ^+_j(w;\hV)\\
&&=E^-(\alpha_{i},z)e_{i}(z)E^+(\al_i,z)E^-(\alpha_{j},w)e_{j}(w)E^+(\al_j,w)\\
&&=\frac{(q^{-b_{ij}+2k}w/z;q^{2k})_\infty(q^{-b_{ij}}w/z;p^*)_\infty}{(q^{b_{ij}+2k}w/z;q^{2k})_\infty
(q^{b_{ij}}w/z;p^*)_\infty}E^-(\alpha_{i},z)e_{i}(z)E^-(\al_j,w)E^+(\al_i,z)
e_{j}(w)E^+(\al_j,w)\\
&&=\frac{(q^{b_{ij}+2k}w/z;q^{2k})_\infty(q^{b_{ij}}w/z;p^*)_\infty}{(q^{-b_{ij}+2k}w/z;q^{2k})_\infty
(q^{-b_{ij}}w/z;p^*)_\infty}E^-(\alpha_{i},z)E^-(\al_j,w)e_{i}(z)e_{j}(w)E^+(\al_i,z)
E^+(\al_j,w)\\
&&=-\frac{w}{z(1-q^{-b_{ij}}w/z)}\frac{(q^{b_{ij}+2k}w/z;q^{2k})_\infty(q^{b_{ij}}z/w;p^*)_\infty}{(q^{-b_{ij}+2k}w/z;q^{2k})_\infty
(p^*q^{-b_{ij}}z/w;p^*)_\infty}\\
&&\qquad\qquad\times E^-(\alpha_{i},z)E^-(\al_j,w)e_{j}(w)e_{i}(z)E^+(\al_i,z)
E^+(\al_j,w)\\
&&=-\frac{w(1-q^{-b_{ij}}z/w)}{z(1-q^{-b_{ij}}w/z)}\frac{(q^{-b_{ij}+2k}z/w;q^{2k})^2_\infty(q^{b_{ij}+2k}w/z;q^{2k})_\infty(q^{-b_{ij}}z/w;p^*)_\infty}{(q^{b_{ij}+2k}z/w;q^{2k})^2_\infty(q^{-b_{ij}+2k}w/z;q^{2k})_\infty
(q^{b_{ij}}z/w;p^*)_\infty}\\
&&\qquad\qquad\times E^-(\al_j,w)e_{j}(w)E^-(\alpha_{i},z)E^+(\al_j,w)e_{i}(z)E^+(\al_i,z)
\\
&&=-\frac{w}{z}\frac{(q^{-b_{ij}}z/w;q^{2k})_\infty(q^{b_{ij}+2k}w/z;q^{2k})_\infty}{(q^{b_{ij}+2k}z/w;q^{2k})_\infty(q^{-b_{ij}}w/z;q^{2k})_\infty}\cZ^+_j(w;\hV)\cZ^+_i(z;\hV).
\en

We also derive \eqref{z2} as follows. 
\be
&&\frac{(q^{b_{ij}+k}w/z;q^{2k})_\infty}{(q^{-b_{ij}+k}w/z;q^{2k})_\infty}\cZ_i^+(z;\hV)\cZ^-_j(w;\hV)\\
&&=\frac{(q^{b_{ij}+k}w/z;q^{2k})_\infty}{(q^{-b_{ij}+k}w/z;q^{2k})_\infty}
E^-(\al_i,z)e_i(z)E^+(\al_i,z)E^-(\al'_j,w)f_j(w)E^+(\al'_j,w)\\
&&=E^-(\al_i,z)E^-(\al'_j,w)e_i(z)f_j(w)E^+(\al_i,z)E^+(\al'_j,w)\\
&&=E^-(\al_i,z)E^-(\al'_j,w)\left[f_j(w)e_i(z)+\frac{\delta_{ij}}{q_i-q^{-1}_i}\left(
\delta\left(q^{-k}\frac{z}{w}\right)\psi^-_i(q^{k/2}w)-\delta\left(q^{k}\frac{z}{w}\right)\psi^+_i(q^{-k/2}w)
\right)\right]\\
&&\qquad\qquad \times E^+(\al_i,z)E^+(\al'_j,w).
\en
Then use 
\be
&&\psi^\pm_i(q^{\mp k/2}w)=K_i^{\pm}E^-(\al_i,q^{\mp k}w)^{-1}E^-(\al'_i,w)^{-1}E^+(\al_i,q^{\mp k}w)^{-1}E^+(\al'_i,w)^{-1}
\en
and the property of the delta function.

To prove the Serre relation (\ref{SerreZ}) for $\cZ^+_j(z)$ we use   
 \eqref{zpe} and \eqref{bosonZ} and obtain  
\bea
&&{e}_{i}(z)=E(\al_i,z)\cZ^{+}_{i}(z;\hV) \lb{eEZ}
\ena 
where we set
\be
E(\al_i,z)&=&E^-(\al_i,z)^{-1}E^+(\al_i,z)^{-1}.
\en 
From \eqref{EpEp}, we have 
\bea
&&E(\al_i,z)E(\al_j,w)\nn\\
&&=\frac{(q^{-2}w/z;q^{2k})_\infty(q^{2}z/w;q^{2k})_\infty}{
(q^{2}w/z;q^{2k})_\infty(q^{-2}z/w;q^{2k})_\infty}
\frac{(p^*q^{-2}w/z;p^*)_\infty(p^*q^{2}z/w;p^*)_\infty}{
(p^*q^{2}w/z;p^*)_\infty(p^*q^{-2}z/w;p^*)_\infty}E(\al_j,w)E(\al_i,z).\nn\\
&&\lb{exchEE}
\ena
Next note that (\ref{serree}) is equivalent to  
\be
0&=&\prod_{1\leq m< l \leq a} 
\frac{(p^*q^{2}{z_{l}}/{z_{m}}; p^*)_{\infty}}
{(p^*q^{-2}{z_{l}}/{z_{m}}; p^*)_{\infty}}
\prod_{1\leq  i \leq a}
\frac{(p^*q^{b_{ij}}{z_{i}}/w; p^*)_{\infty}}
{(p^*q^{-b_{ij}}{z_{i}}/w; p^*)_{\infty}}\\
&&\times
\sum_{\sigma\in S_{a}}\prod_{1\leq m< l \leq a \atop 
\sigma^{-1}(m)>\sigma^{-1}(k)} 
\frac{(p^*q^{-2}{z_{l}}/{z_{m}}; p^*)_{\infty}}
{(p^*q^{2}{z_{l}}/{z_{m}}; p^*)_{\infty}}
\frac{(p^*q^{2}{z_{\sigma(l)}}/{z_{\sigma(m)}}; p^*)_{\infty}}
{(p^*q^{-2}{z_{\sigma(l)}}/{z_{\sigma(m)}}; p^*)_{\infty}}\nn\\
&\times&\sum_{s=0}^{a}(-1)^{s}\left[\begin{array}{c}
a\nn\\
s\cr
\end{array}\right]_{i}\prod_{1\leq  m \leq s}
\frac{(p^*q^{b_{ij}}{w}/{z_{\sigma(m)}}; p^*)_{\infty}}
{(p^*q^{-b_{ij}}{w}/{z_{\sigma(m)}}; p^*)_{\infty}}
\frac{(p^*q^{-b_{ij}}{z_{\sigma(m)}}/{w}; p^*)_{\infty}}
{(p^*q^{b_{ij}}{z_{\sigma(m)}}/{w}; p^*)_{\infty}}\\
&&\times e_{i}(z_{\sigma(1)})\cdots e_{i}(z_{\sigma(s)})e_{j}(w)e_{i}(z_{\sigma(s+1)})\cdots e_{i}(z_{\sigma(a)}).
\label{eq1}
\en
Substitute \eqref{eEZ}  into this, 
and move all $E(\al_i,z_j)$ and $E(\al_j,w)$ to the left. Then we get
\bea
0&=&
\prod_{1\leq m< l \leq a} 
\frac{(p^*q^{2}{z_{l}}/{z_{m}}; p^*)_{\infty}}
{(p^*q^{-2}{z_{l}}/{z_{m}}; p^*)_{\infty}}
\prod_{1\leq  i \leq a}
\frac{(p^*q^{b_{ij}}{z_{i}}/w; p^*)_{\infty}}
{(p^*q^{-b_{ij}}{z_{i}}/w; p^*)_{\infty}}\nn\\
&&\times
\sum_{\sigma\in S_{a}}\prod_{1\leq m< l \leq a \atop 
\sigma^{-1}(m)>\sigma^{-1}(k)} 
\frac{(p^*q^{-2}{z_{l}}/{z_{m}}; p^*)_{\infty}}
{(p^*q^{2}{z_{l}}/{z_{m}}; p^*)_{\infty}}
\frac{(p^*q^{2}{z_{\sigma(l)}}/{z_{\sigma(m)}}; p^*)_{\infty}}
{(p^*q^{-2}{z_{\sigma(l)}}/{z_{\sigma(m)}}; p^*)_{\infty}}\nn\\
&\times&\sum_{s=0}^{a}(-1)^{s}\left[\begin{array}{c}
a\nn\\
s\cr
\end{array}\right]_{i}\prod_{1\leq  m \leq s}
\frac{(p^*q^{b_{ij}}{w}/{z_{\sigma(m)}}; p^*)_{\infty}}
{(p^*q^{-b_{ij}}{w}/{z_{\sigma(m)}}; p^*)_{\infty}}
\frac{(p^*q^{-b_{ij}}{z_{\sigma(m)}}/{w}; p^*)_{\infty}}
{(p^*q^{b_{ij}}{z_{\sigma(m)}}/{w}; p^*)_{\infty}}\\
&&\times\mathfrak{\varepsilon}(z_{\sigma(1)},\cdots,z_{\sigma(s)},w,z_{\sigma(s+1)},\cdots,z_{\sigma(a)})\nn\\
&&\times \cZ^+_{i}(z_{\sigma(1)};\hV)\cdots \cZ^+_{i}(z_{\sigma(s)};\hV)\cZ^+_{j}(w;\hV)\cZ^+_{i}(z_{\sigma(s+1)};\hV)\cdots \cZ^+_{i}(z_{\sigma(a)};\hV),
\label{eq1}
\ena
where we set 
\be
&&\mathfrak{\varepsilon}
(z_{\sigma(1)},\cdots,z_{\sigma(s)},w,z_{\sigma(s+1)},\cdots,z_{\sigma(a)})\nn\\
&&=E(\al_i,z_{\sigma(1)})\cdots E(\al_i,z_{\sigma(s)}) E(\al_j,w) E(\al_i,z_{\sigma(s+1)}) \cdots
E(\al_i,z_{\sigma(a)}).\lb{productbosons}
\en
Then moving $E(\al_j,w)$ to the left end by 
 (\ref{exchEE}), we have 
\be
&&\mathfrak{\varepsilon}(z_{\sigma(1)},\cdots,z_{\sigma(s)},w,z_{\sigma(s+1)},\cdots,z_{\sigma(N)})\nn\\
&&=\prod_{1\leq  i \leq s}
\frac{(q^{-b_{ij}}w/z_{\sigma(i)};q^{2k})_{\infty}
(q^{b_{ij}}z_{\sigma(i)}/w;q^{2k})_{\infty}}{(q^{b_{ij}}w/z_{\sigma(i)};q^{2k})_{\infty}(q^{-b_{ij}}z_{\sigma(i)}/w;q^{2k})_{\infty}}
\frac{(p^*q^{-b_{ij}}w/z_{\sigma(i)};p^\ast)_{\infty}(p^*q^{b_{ij}}z_{\sigma(i)}/w;p^\ast)_{\infty}}
{(p^*q^{b_{ij}}w/z_{\sigma(i)};p^\ast)_{\infty}(p^*q^{-b_{ij}}z_{\sigma(i)}/w;p^\ast)_{\infty}}\nn\\
&&\qquad\qquad\times\mathfrak{\varepsilon}
(w,z_{\sigma(1)},\cdots,z_{\sigma(a)}).\label{vep1}
\en
Substituting this  into (\ref{eq1}), we can factor out
 $\mathfrak{\varepsilon}(w,z_{\sigma(1)},\cdots,z_{\sigma(a)})$
 from $\sum_{s=0}^{a}$. 
Then exchanging the order of $E(\al_i,z_l)$'s by \eqref{exchEE}, we have 
\be
&&\mathfrak{\varepsilon}
(w,z_{\sigma(1)},\cdots,z_{\sigma(a)})\nn\\
&&=\prod_{1\leq m< l \leq a\atop
\sigma^{-1}(m)>\sigma^{-1}(k)}
\frac{(q^{-2}{z_{l}}/{z_{m}};q^{2k})_{\infty}(q^{2}{z_{\sigma(l)}}/{z_{\sigma(m)}};q^{2k})_{\infty}}
{(q^{2}{z_{l}}/{z_{m}};q^{2k})_{\infty}(q^{-2}{z_{\sigma(l)}}/{z_{\sigma(m)}};q^{2k})_{\infty}}
\frac{(p^*q^{2}{z_{l}}/{z_{m}};p^\ast)_{\infty}
(q^{-2}p^{\ast}{z_{\sigma(l)}}/{z_{\sigma(m)}};p^\ast)_{\infty}}
{(p^*q^{-2}{z_{l}}/{z_{m}};p^\ast)_{\infty}
(p^*q^{2}{z_{\sigma(l)}}/{z_{\sigma(m)}};p^\ast)_{\infty}}
\nn\\
&&\qquad\qquad\times \mathfrak{\varepsilon}(w,z_{1},\cdots,z_{a}).
\label{e2}
\en
Substituting this into (\ref{eq1}), we can factor out $\mathfrak{\varepsilon}
(w,z_{1},\cdots,z_{a})$ 
completly from $\sum_{\sigma\in S_{a}}$. 
Multiply 
\be
&&\prod_{1\leq m< l \leq a}
\frac{(q^{2}{z_{l}}/{z_{m}};q^{2k})_{\infty}}{(q^{-2}{z_{l}}/{z_{m}};q^{2k})_{\infty}}
\prod_{1\leq m\leq a}\frac{(q^{-b_{ij}}{z_{m}}/w;q^{2k})_{\infty}}{(q^{b_{ij}}{z_{m}}/w;q^{2k})_{\infty}},
\en
 and drop the overall factor depending on $p^*$,  one  gets the desired relation.
One can prove the $\cZ^-_j(z;\hV)$ case in the same way.

\begin{dfn}\label{z-algebra}
For $k\in \C^\times$ and $(\hV,\pi)\in \fC_k$, we call 
the $H$-subalgebra of $\cD_{H,\hV}$
 generated by $\cZ_{i,m}^\pm(\hV)$, $K^\pm_i\ (i\in I, m\in \Z), \cM_{H^*}$ and $\hd$ 
the dynamical quantum $Z$-algebra ${\cZ}_{\hV}$ 
 associated with  $(\hV,\pi)$. 
\end{dfn}

\subsection{The universal algebra $\hcZ_{k}$}
Using the relations in Theorem \ref{zz}, we define the universal dynamical quantum $Z$-algebra as follows. 

\begin{dfn}\label{univz-algebra}
Let $\hcZ^\pm_{i,m}
\ (i\in I, m\in \Z)$ be abstract symbols.  
We set $\hcZ^\pm_i(z)=\sum_{m\in \Z}\hcZ^\pm_{i,m}z^{-m}$. 
We define  the universal dynamical quantum $Z$-algebra $\hcZ_{k}$ to be a topological algebra over $\K[[q^{2{k}}]]$ generated 
by $\hcZ^\pm_{i,m}, K^\pm_i \ (i\in I, m\in \Z), \hd, \cM_{H^*}$
subject to the relations obtained by replacing 
$\cZ^\pm_i(z;\hV)$ by $\hcZ^\pm_i(z)$  in Theorem \ref{zz}. 
\end{dfn}
We treat the relations as formal Laurent series in $z, w$ and  $z_j$'s in a similar way to those of $U_{q,p}(\gh)$ 
in Sec.2.1. The defining relations are well-defined in the $q^{2k}$-adic topology.

\begin{prop}
${\cZ}_{k}$ is an $H$-algebra with the same $\mu_l,\mu_r$ as in $U_{q,p}(\gh)$. 
\end{prop}

Note that for $(\hV,\hpi)\in \sC_k$ we extend $\hpi$ to the map $\hpi: {\cZ}_k\to \cD_{H,\hV}$ by 
$\hpi(\hcZ^\pm_{i,m})=\cZ^\pm_{i,m}(\hV) $. Then $\hV$ is a ${\cZ}_k$-module by $\hpi$.

\begin{dfn}
For $k\in \C^\times$, we denote by $\sD_k$ the full subcategory of the category of 
${\cZ}_k$-modules consisting of those modules $(\cW,\sigma)$ such that
\begin{itemize}
\item[(i)] $\cW$ has level $k$.
\item[(ii)] $\cW=\bigsqcup_{\omega\in \C} \cW_\omega$, where  $\cW_\omega=\{w\in \cW\ |\ \sigma(\hd) w=\omega w\ \}$
\item[(iii)] For every $\omega \in \C$, there exists $n_0\in \N$ such that for all 
$n>n_0$, $\cW_{\omega+n}=0$.   
\end{itemize}
\end{dfn}

Let us consider $(\hV,\hpi)\in \sC_k$.   
Following Lepowsky and Wilson\cite{LW}, 
we define the vacuum space $\Omega_{\hV}$ by
\be
&&\Omega_{\hV}=\{ v\in \hV\ |\ \hpi(\al_{i,n}) v=0\quad \forall i\in I,\ n\in \Z_{>0}\ \}. 
\en
From Theorem \ref{zz}, $\Omega_{\hV}$ is stable under the 
action of $\hcZ_{\hV}$.  For a morphism $f:\hV\to \hV'$ in $\sC_k$, we have 
\be
&&f(\Omega_{\hV})\subset \Omega_{\hV'}. 
\en

\begin{prop}
For $(\hV,\hpi)\in \sC_k$, there is a unique representation $\sigma$ of ${\cZ}_k$ on $\Omega_{\hV}$ 
such that $(\Omega_{\hV},\sigma)\in \sD_k$, 
\be
&&\sigma(K^\pm_i)=\hpi(K^\pm_i),\quad \sigma(\cZ^\pm_{i,m})=\cZ^\pm_{i,m}(\hV)\quad \forall i\in I, m\in \Z.
\en
\end{prop}
We hence  
 define a functor $\Omega: \sC_k\to \cD_k$ by
\be
&&\Omega(\hV,\hpi)=(\Omega_{\hV},\sigma),\quad \Omega({f})=f|_{\Omega_{\hV}}:\Omega_{\hV}\to\Omega_{\hV'}.
\en

\subsection{The functor $\Lambda$}
We define a reverse functor $\Lambda:\fD_k\to \fC_k$ as follows. 
Let $(\cW,\hs)\in \fD_k$ be a ${\cZ}_k$-module. We define $U_{q,p}(\cH)$-module 
$\Ind\, \cW$ by requiring $\al_{i,m}\cdot \cW=0$ and 
\be
\Ind \, \cW=U_{q,p}(\cH)\otimes_{U_{q,p}(\cH_+)}\cW.
\en
Let $\F_{\al,k}$ be the level-$k$ Fock module defined in Sec.3.1.  We have a natural isomorphism $ 
\F_{\al,k}\otimes_\C \cW \cong\Ind\ \cW$ by $(u\otimes 1_k)\otimes w\mapsto u\otimes w$\cite{LW}. We thus identify the $U_{q,p}(\cH)$-module $\Ind\, \cW $ with $\F_{\al,k}\otimes_\C \cW$,  with the action $\hpi$ of ${U}_{q,p}(\cH)$ 
\be
&&\hpi(c)=1\otimes c,
\quad \hpi(K^\pm_i)=1\otimes \sigma(K^\pm_i),
\quad \hpi(\al_{i,m})=\al_{i,m}\otimes 1.   
\en
For $(\cW,\sigma)\in \fD_k$ and $\Ind\, \cW=\F_{\al,k}\otimes_\C \cW$, 
we define $e'_j(z), f'_j(z)\in \cD_{H, \Ind\, \cW}[[z,z^{-1}]]$ by
\be
&&e'_j(z)=E^-(\al_j,z)^{-1}E^+(\al_j,z)^{-1}\otimes \sigma(\cZ^{+}_{j}(z)) ,\\
&&f'_j(z)=E^-(\al'_j,z)^{-1}E^+(\al'_j,z)^{-1}\otimes \sigma(\cZ^{-}_{j}(z)).
\en
These are well-defined elements of $\cD_{H,\Ind\, \cW}[[z,z^{-1}]]$. 
By a similar argument to the proof of Theorem \ref{zz} one can show that $e'_j(z)$ and $f'_j(z)$ satisfy the defining relations of $U_{q,p}(\gh)$ with $c=k$. 
We hence extend $\hpi: U_{q,p}(\cH)\to \cD_{H,\Ind\, \cW}$ to 
$\hpi: U_{q,p}(\gh)\to \cD_{H,\Ind\, \cW}$ as an $H$-algebra homomorphism by 
\be
&&\hpi(e_j(z))=e'_j(z),\quad  \hpi(f_j(z))=f'_j(z),\\
&&\hpi(\hd)=\hd\otimes 1+1\otimes \sigma(\hd). 
\en
By construction, the latter map is uniquely determined. 
\begin{prop}\lb{IndW}
For $(\cW,\sigma)\in \fD_k$, there is a unique level-$k$ $U_{q,p}(\gh)$-module 
$(\Ind\, \cW, \hpi)\in \fC_k$.  
\end{prop}

We thus reach the following definition. 
\begin{dfn}
We define a functor $\Lambda:\fD_k\to \fC_k$ by
\begin{itemize}
\item[(i)]\ $\Lambda(\cW,\sigma)=(\Ind\, \cW,\hpi)$
\item[(ii)]\ For a morphism $f:\cW\to \cW'$ in $\fD_k$, define $\Lambda(f):\Ind\, \cW\to \Ind\, \cW'$ 
to be the induced $U_{q,p}(\cH)$-module map. Then $\Lambda(f)$ is a $U_{q,p}(\gh)$-module map. 
\end{itemize}
\end{dfn}

We obtain the following theorem analogously to the case of the affine Lie algebras\cite{LW}. 
\begin{thm}\lb{IrrIndW}
For $k\in\C^\times$, the two categories $\fC_k$ and $\fD_k$ are equivalent by 
the functors $\Omega:\fC_k\to \fD_k$ and $\Lambda:\fD_k\to \fC_k$. 
In particular, the level-$k$ $U_{q,p}(\gh)$-module $\Ind\, \cW=\F_{\al,k}\otimes_\C \cW \in \fC_k$ is irreducible 
if and only if $\cW\in \fD_k$ is an irreducible ${\cZ}_k$-module.   
\end{thm}

\section{The Induced $U_{q,p}(\gh)$-Modules }
In this section we give a simple realization of the dynamical quantum $Z$-algebra 
$\cZ_k$ in terms of the quantum $Z$-algebra $Z_k$ associated with $U_q(\gh)$ and  
 construct the level-$k$ induced $U_{q,p}(\gh)$-modules. 
We also give some examples of the level-1 irreducible representations.

\subsection{The quantum $Z$-algebra ${Z}_k$ associated with $U_q(\gh)$}
One can apply the arguments similar to those in Secs.3.1-3.3 to the quantum affine algebra $U_q(\gh)$ in 
the Drinfeld realization and define the corresponding  
quantum $Z$-algebras ${Z}_{V}$  
associated with the level-$k$ $U_q(\gh)$-module $V$\cite{Jing} and the universal one 
${Z}_k$. See Appendix \ref{Uqgh}. We also denote by $C_k$ and $D_k$ 
the $U_q(\gh)$ counterparts of the categories $\fC_k$ and $\fD_k$.

Comparing the defining relations of $\cZ_{k}$ with those of ${Z}_{k}$, 
we obtain the following isomorphism.    
\begin{prop}\lb{cZHga}
We have the isomorphism 
\be
&&{\cZ}_{k}\cong (\K\otimes_\C {Z}_{k})\sharp\C[\cR_Q]
\en
as an $H$-algebra by  
\be
&&  \hcZ_{j,m}^+ \mapsto  {Z}_{j,m}^+e^{-Q_{\al_j}},\quad  \hcZ_{j,m}^- \mapsto  {Z}_{j,m}^-,\\
&&K^\pm_i\mapsto q_i^{\mp h_i}e^{-Q_{\al_j}}\  (i\in I, m\in \Z),\quad \hd \mapsto \bar{d},  
\en
where $Z^\pm_{j,m}$ denotes the generators in $Z_k$ (Definition \ref{z-algebra}). 
\end{prop}

\begin{thm}\lb{cWHVQ}
For $(W, \bar{\sigma})\in D_k$ and generic $\mu\in \h^*$, there is a dynamical representation ${\hs}$ of ${\cZ}_k$ 
on ${\cW_{H,Q}}(\mu):=(\K\otimes_\C W)\otimes_{\C} e^{Q_{\bar{\mu}}}\C[\cR_Q]$  such that $(\cW_{H,Q}(\mu),\hs) \in \fC_k$ and 
\be
&&\hs({\hcZ^+_{j,m}})=\bs({Z^+_{j,m}})\otimes e^{-Q_{\al_j}}, \quad \hs({\hcZ^-_{j,m}})=\bs({Z^-_{j,m}})\otimes 1,\\
&&\hs(K^\pm_j)=\bs(q_j^{\mp h_j})\otimes e^{-Q_{\al_j}},\quad \hs(\hd)=\bs(\bd)\otimes 1+1\otimes P_d, 
\en
where $P_d$ denotes a $\C$-linear operator on $1\otimes e^{Q_{\bar{\mu}}}\C[\cR_Q]$ such that 
\be
&&[1\otimes P_d, \hs({\hcZ^\pm_{j,m}})]=0. 
\en
\end{thm}

\begin{prop}\lb{IrrW}
The representation $(\cW_{H,Q}(\mu),\hs)$ of 
 $\hcZ_k$ is irreducible  if and only if $W$ is an irreducible $Z_k$-module.  
\end{prop}

From this and Theorem \ref{IrrIndW}, we obtain: 
\begin{prop}\lb{IrrStandard}
For a $Z_k$-module $(W,\bs)\in D_k$  and generic $\mu\in \h^*$, let $(\cW_{H,Q}(\mu)
,\s)$ be the $\cZ_k$-module constructed in Theorem \ref{cWHVQ} and 
$\Ind\, \cW_{H,Q}(\mu)=\F_{\al,k}\otimes_\C \cW_{H,Q}(\mu)$ be the level-$k$ induced $U_{q,p}(\gh)$-module given in  Proposition \ref{IndW}. 
Then  $(\Ind\, \cW_{H,Q}(\mu), \pi )$ is irreducible if and only if $(W,\bs)$ is  irreducible. 
 \end{prop}

\subsection{Examples of the irreducible representations}\lb{level1Ex}

We here give some examples of the level-1 irreducible induced representations of $U_{q,p}(\gh)$  
of   types $\gh=A_l^{(1)}, D_l^{(1)}, E_6^{(1)}, E_7^{(1)}, E_8^{(1)}$ and $B_l^{(1)}$.

\noindent
\subsubsection{The simply laced case :}\lb{ADE} 

Let $\C[\cQ]$ be the group algebra of the root lattice $\cQ=\oplus_i \Z\al_i$ 
with the central extension:  
\be
&&e^{ \alpha_i}e^{ \alpha_j}=(-1)^{ (\alpha_i, \alpha_j)}e^{ \alpha_j}e^{ \alpha_i}\qquad 
(i,j\in I). 
\en
Let us consider the fundamental weight $\Lambda_a $ of $\gh$ with $0\leq a\leq l$ for $A^{(1)}_l$, $a=0,1,l-1,l$  for $D_l^{(1)}$, $a=0,1,2$ for $E_6^{(1)}$, $a=0,1$ for $E_7^{(1)}$, $a=0$ for $E_8^{(1)}$. 
\begin{thm}\cite{FJ,Jing}
An  inequivalent set of the level-1 irreducible $Z_1(\gh)$-{modules} is given by 
$W(\Lambda_a)=e^{\bar{\Lambda}_a}\C[\cQ]$, on which 
  the actions of $Z_j^\pm(z)$ 
 are given by 
\bea
&&Z_{j}^{\pm}(z)=e^{\pm \alpha_j}z^{\pm h_{j}+1} \lb{qZop}
\ena
with  
\be
z^{\pm h_{i}}e^{\pm \alpha_j}e^{\bar\Lambda_a}=z^{\pm(\al_{i}^\vee,\alpha_j+\bar{\Lambda}_a)}e^{\pm \alpha_j}
e^{\bar\Lambda_a} \qquad (i,j\in I).
\en
\end{thm}

Then for generic  $\mu\in \h^*$,    
we have from Theorem \ref{cWHVQ} a level-1 irreducible $\cZ_1(\gh)$ module 
 $\cW_{H,Q}(\Lambda_a,\mu):=(\K\otimes_\C W(\Lambda_a))\otimes e^{Q_{\bar{\mu}}}\C[\cR_Q]$ 
with the action given by 
\bea
&&\cZ_{j}^{+}(z)=Z_{j}^{+}(z)\otimes e^{-Q_{\al_j}}
,\quad \cZ_{j}^{-}(z)=Z_{j}^{-}(z)\otimes 1 \lb{qcZop}
. 
\ena
Then from Proposition \ref{IrrStandard} we obtain:
\begin{thm}
 A level-1 irreducible highest weight representations of $U_{q,p}(\gh)$ 
is given by $\hV(\Lambda_a+\mu,\mu):=\Ind\, {\cW}_{H,Q}(\Lambda_a,\mu)$ with the highest weight $(\Lambda_a+\mu, \mu)$:  
\be
&&\hV(\Lambda_a+\mu,\mu)=\F_{\al,1}\otimes\cW_{H,Q}(\Lambda_a,\mu)
=\bigoplus_{\gamma,\kappa\in \cQ}\F_{\gamma,\kappa}(\Lambda_a,\mu), 
\en
where
\be
&&\F_{\gamma,\kappa}(\Lambda_a,\mu)=\K\otimes_\C  (\F_{\al,1}\otimes e^{\bar{\Lambda}_a+\gamma}) \otimes e^{Q_{\bar{\mu}+\kappa}}, 
\en
The highest weight vector is $ 1_1\otimes e^{\bar{\Lambda}_a}\otimes e^{Q_{\bar{\mu}}}$. 
The derivation operator $\hd$ is realized as  
\be
\hd&=&-\frac{1}{2}\sum_{j=1}^lh_jh^j-N^\al+\frac{1}{2r^*}\sum_{j=1}^l(P_j+2)P^j-\frac{1}{2r}\sum_{j=1}^l((P+h)_j+2)(P+h)^j,\\
N^\al&=&\sum_{j=1}^l\sum_{m\in\Z_{>0}}\frac{m^2}{[m]}\frac{1-p^{*m}}{1-p^m} q^m\al_{j,-m}A^j_m,
\en
where $r, r^*\in \C^\times$, and  $A^j_m$ are the fundamental weight type elliptic bosons given in Sec.5.1. 
\end{thm}

One can easily calculate the character of $\hV(\Lambda_a,\mu)$: 
\be
ch_{\hV(\Lambda_a+\mu,\mu)}&=&\tr_{\hV(\Lambda_a+\mu,\mu)}q^{-d-\frac{c({W}(\g))}{24}}=\sum_{\gamma,\kappa\in \cQ}ch_{\F_{\gamma,\kappa}(\Lambda_a,\mu)},\\
ch_{\F_{\gamma,\kappa}(\Lambda_a,\mu)}&=&\frac{1}{\eta(q)^l}q^{\frac{1}{2rr^*}|r(\bar{\mu}+\kappa+\bar{\rho})-r^*(\bar{\Lambda}_a+\bar{\mu}+\gamma+\kappa+\bar{\rho})|^2}.
\en
Here $c(W(\g))=l(1-\frac{g(g+1)}{rr^*})$, and $\eta(q)$ denotes Dedekind's $\eta$-function 
given by
\be
&&\eta(q)=q^{\frac{1}{24}}(q;q)_\infty.
\en
One should note that the character $ch_{\F_{\gamma,\kappa}(\Lambda_a,\mu)}$ coincides with 
the one of the Verma module of the $W(\g)$-algebras for $\g=A_l, D_l, E_6, E_7, E_8$ with the highest weight $h=\frac{1}{2rr^*}|r(\bar{\mu}+\kappa+\rho)-r^*(\bar{\Lambda}_a+\bar{\mu}+\gamma+\kappa+\rho)|^2$ and the central charge 
$c(W(\g))$. 
In fact, for $\gh=A^{(1)}_l$ case, for example, one can construct an action of the deformed 
$W(A_l)$ algebra on $\F_{\gamma,\kappa}(\Lambda_a,\mu)$ 
explicitly. 
\begin{thm}\lb{dWA}\cite{FeFr,AKOS}
For $p=q^{2r}$ and $p^*=pq^{-2}=q^{2r^*}$, i.e. $r^*=r-1$, the deformed $W(A_l)$-algebra acts on 
$\F_{\xi,\kappa}(\Lambda_i+\mu,\mu)$ by
\be
&&\Lambda_j(z)=:\exp\left\{ \sum_{m\not=0}(q^m-q^{-m})(1-p^{*m})\cE^{+j}_m(q^jz)^{-m}\right\}:\otimes\ 
p^{*h_{\bar{\bep}_j}}\quad (1\leq j\leq l),\\
&&T_n(z)=\sum_{1\leq j_1<\cdots<j_n\leq l}: \Lambda_{j_1}(z)\Lambda_{j_2}(zq^{-2})\cdots \Lambda_{j_n}(zq^{-2(n-1)}):\quad 
(1\leq n\leq l). 
\en
Here $\cE^{+j}_m$ denotes the orthonormal basis type elliptic boson given in \eqref{ONboson}, and {$:\quad :$} denotes the normal ordering of the enclosed expression such that the operators $\cE^{\pm j}_m$ for $m<0$ are to be placed to the left of the operators $\cE^{\pm j}_m$ for $m>0$.{} 
In addition, the level-1 elliptic currents $e_j(w)$ and $f_j(w)$ of $U_{q,p}(A^{(1)}_l)$ obtained from Proposition \ref{IndW}, \eqref{qZop} and \eqref{qcZop} 
are the screening currents of the deformed $W(A_l)$-algebra, i.e. 
they commute with $T_n(z)$ up to a total difference. 
\end{thm}
See also \cite{JKOS, KojimaKonno,KK04}. A similar statement is valid also for the deformed $W(D_l)$\cite{FrRe} and $U_{q,p}(D^{(1)}_l)$.  
We also expect that for $r\in \Z_{>0}$ satisfying $r>g+1$ and for a level-$(r-g-1)$ 
dominant integral weight $\mu$, the space $\F_{\gamma,\kappa}(\Lambda_a,\mu)$ 
 becomes completely degenerate with respect to the action of the corresponding deformed 
$W(\g)$-algebra \cite{FeFr,AKOS,FrRe}, although the $E_{6,7,8}$-type deformed $W$ algebras have not yet been 
constructed explicitly. In order to get the irreducible module one should make 
the BRST-resolution in terms of the BRST-charge constructed from the half currents of $U_{q,p}(\gh)$. 
An explicit demonstration for the $A_1^{(1)}$ case has been discussed in \cite{KKW}. 

\noindent
{\it Remark.} In Theorem \ref{dWA}, we assumed $p=q^{2r}$ in order to make a connection to the deformed 
$W(A_l)$-algbera. The same relation arises naturally when one considers the finite dimensional representations 
of the universal elliptic dynamical $\cR$ matrices\cite{JKOStg,Konno06}.    
     
\noindent
\subsubsection{The $B_l^{(1)}$ case}  

We follow the work \cite{LP84} and its quantum analogues \cite{Be,JinMis} with a slight modification in the Ramond sector according to \cite{Idz}. 
Let  $e^{\alpha_i}\ (i\in I)$ be the generators of the group algebra $\C[\cQ]$ with the following central extension. 
 \be
&&e^{ \alpha_i}e^{ \alpha_j}=(-1)^{(\alpha_i, \alpha_j)+(\alpha_i, \alpha_i)(\alpha_j, \alpha_j)} e^{\alpha_j}e^{ \alpha_i}
\en
As before we regard $h_i \ (i\in I)$ as an operator such that 
\be
&& z^{\pm h_{i}}e^{ \alpha_j}=z^{\pm(\alpha^\vee_{i},\alpha_{j})}e^{ \alpha_j}z^{\pm h_{i}}
\en
We also need the  Neveu-Schwartz ($NS$) fermion
$\{\Psi_n | n \in {\mathbb Z}+\frac{1}{2}\}$
and the Ramond ($R$) fermion
$\{\Psi_n | n \in {\mathbb Z}\}$ 
satisfying the following anti-commutation relations.
\begin{eqnarray*}
~\{\Psi_m,\Psi_n\}=\delta_{m+n,0}\cN({q^{m}+q^{-m}})
\end{eqnarray*}
with $\cN=1/(q^{\frac{1}{2}}+q^{-\frac{1}{2}})$. 
We  define 
\be
&&{\cal F}^{NS}={\mathbb C}[\Psi_{-\frac{1}{2}},\Psi_{-\frac{3}{2}} 
\cdots],\qquad 
\widetilde{\F}^{R}=\C[\Psi_{-1},\Psi_{-2},...]
\en
and their submodules $\F^{NS, R}_{even}$ (reps. $\F^{NS, R}_{odd}$)  
generated by the even (reps. odd) number of $\Psi_{-m}$'s.  
One should note that for the $R$ fermion 
$\Psi_0^2=\cN$ and $\{\Psi_m,\Psi_0\}=0$ for $m\not=0$.   
So we have two degenerate vacuum states $1$ and $\Psi_01$. 
We hence consider the extended space 
\be
&&\widehat{\F}^R=\widetilde{\F}^R\otimes \C^2
\en
and realize the  $R$-fermions by   
\be
&&\widehat{\Psi}_m=\Psi_m\otimes \left(\mmatrix{1&0\cr 0&-1\cr}\right) 
\quad (m\in \Z_{\not =0} )
,\qquad
\widehat{\Psi}_0=\cN^{\frac{1}{2}}(1\otimes \left(\mmatrix{0&1\cr 1&0\cr}\right)).  
\en
Note that $\{\hPsi_m,\hPsi_n\}=\delta_{m+n,0}\cN({q^{m}+q^{-m}})$. We set 
\be
&&\F^R=\F^R_{even}
\otimes \C\left(\mmatrix{1\cr 1\cr}\right)
\oplus \F^R_{odd}\otimes\C\left(\mmatrix{1\cr -1\cr}\right).\lb{RFock}
\en
The action of $\Psi_{m}$ on $\F^{NS}$ is given by
\be
&&\Psi_{-m}\cdot u= \Psi_{-m}u,\qquad 
\Psi_{m}\cdot u=\{\Psi_{m}, u\}\qquad (m\in\Z_{>0}),
\en
where $u\in\F^{NS}$, whereas $\widehat{\Psi}_{m}$ acts on $\F^{R}$ as
\be
&&\widehat{\Psi}_{-m}\cdot u\otimes v= \Psi_{-m}u\otimes \left(\mmatrix{1&0\cr 0&-1\cr}\right)v\ (m\in\Z_{>0}),
\quad \widehat{\Psi}_{0}\cdot u\otimes v= u\otimes \left(\mmatrix{0&1\cr 1&0\cr}\right)v,\nn\\
&&\widehat{\Psi}_{m}\cdot u\otimes v=\{\Psi_{m}, u\}\otimes \left(\mmatrix{1&0\cr 0&-1\cr}\right)v\qquad (m\in\Z_{>0}),
\en 
where $u\in\widetilde{\F}^{R},\ v\in \C^2$. 

Let us define the fermion fields    
$\Psi^{NS}(z)$ and $\Psi^{R}(z)$ by
\begin{eqnarray*}
&&\Psi^{NS}(z)=\sum_{n \in {\mathbb{Z}}+\frac{1}{2}}
\Psi_n z^{-n},\qquad
\Psi^{R}(z)=\sum_{n \in {\mathbb{Z}}}
\widehat{\Psi}_n z^{-n}.
\end{eqnarray*}
One can derive the following operator product expansions. 
\be
\Psi(z)\Psi(w)=:\Psi(z)\Psi(w):+<\Psi(z)\Psi(w)>,
\en
where
\be
<\Psi(z)\Psi(w)>=\left\{
\mmatrix{
    \frac{(zw)^{1/2}(z-w)}{(z-qw)(z-q^{-1}w)} & \mbox{for NS}\cr 
    \cN\frac{(z-w)(z+w)}{(z-qw)(z-q^{-1}w)} & \mbox{for R}.\cr
}
\right.
\en
Then the quantum $Z$-algebra $Z_1(B^{(1)}_l)$ is realized as follows\cite{Jing}.  
\be
&&Z_i^{\pm}(z)=e^{\pm \alpha_i}z^{\pm h_{i}+1}\qquad (1\leq i\leq l-1),\nn\\
&&Z_l^{\pm}(z)=\frac{1}{\cN^{1/2}}\Psi(z)e^{\pm \alpha_l}z^{\pm d_lh_{l}+d_l}. 
\en

There are three irreducible $Z_1(B^{(1)}_l)$-modules given by 
\be
&&W(\Lambda_0)= \F^{NS}_{even}\otimes\C[\cQ_0]\oplus \F^{NS}_{odd}\otimes \C[\cQ_0]e^{\bar{\Lambda}_1},\\
&&W(\Lambda_1)=\F^{NS}_{even}\otimes \C[\cQ_0]e^{\bar{\Lambda}_1}\oplus \F^{NS}_{odd}\otimes \C[\cQ_0],\\
&&W(\Lambda_l)=\F^R\otimes \C[\cQ]e^{\bar{\Lambda}_l}
\cong \F^R\otimes \C[\cQ_0]e^{\bar{\Lambda}_l}\oplus \F^R\otimes \C[\cQ_0] e^{\bar{\Lambda}_1+\bar{\Lambda}_l},
\en
where $\cQ_0$ denotes the sublattice of $\cQ$ generated by the long roots. 
For generic $\mu\in \h^*$ and $a=0,1,l$, we set $\cW_{H,Q}(\Lambda_a,\mu)=(\K\otimes_\C W(\Lambda_a))\otimes e^{Q_{\bar{\mu}}}\C[\cR_Q]$. 
From Proposition \ref{IndW} we have the following three level-1 irreducible $U_{q,p}(\hat{B}_{l}^{(1)})$-modules with the higest weight $(\Lambda_a+\mu,\mu)$: 
\be
\hV(\Lambda_{a}+\mu,\mu)&=&\F_{\al,1}\otimes_\C \cW_{H,Q}(\Lambda_a,\mu)
\\
&=&\bigoplus_{\gamma\in \cQ_0,\kappa\in \cQ}\bigoplus_{\la\in {max}(\Lambda_a)\atop mod\ Q_0+\C\delta} 
\F_{\la,\gamma,\kappa}(\Lambda_a,\mu),
\en
where 
\be
&&\F_{\Lambda_0,\gamma,\kappa}(\Lambda_0,\mu)=\K\otimes_\C(\F_{\al,1}\otimes \F^{NS}_{even}\otimes e^\gamma)\otimes e^{Q_{\bar{\mu}+\kappa}},\\
&&\F_{\Lambda_1,\gamma,\kappa}(\Lambda_0,\mu)=\K\otimes_\C(\F_{\al,1}\otimes\F^{NS}_{odd}\otimes e^{\bar{\Lambda}_1+\gamma})\otimes e^{Q_{\bar{\mu}+\kappa}},\\
&&\F_{\Lambda_1,\gamma,\kappa}(\Lambda_1,\mu)=\K\otimes_\C(\F_{\al,1}\otimes\F^{NS}_{even}\otimes e^{\bar{\Lambda}_1+\gamma})\otimes e^{Q_{\bar{\mu}+\kappa}},\\
&&\F_{\Lambda_0,\gamma,\kappa}(\Lambda_1,\mu)=\K\otimes_\C(\F_{\al,1}\otimes\F^{NS}_{odd}\otimes e^{\gamma})\otimes e^{Q_{\bar{\mu}+\kappa}},\\
&&\F_{\Lambda_l,\gamma,\kappa}(\Lambda_l,\mu)=\K\otimes_\C(\F_{\al,1}\otimes\F^{R}\otimes e^{\bar{\Lambda}_l+\gamma})\otimes e^{Q_{\bar{\mu}+\kappa}},\\
&&\F_{\Lambda_l-\al_l,\gamma,\kappa}(\Lambda_1,\mu)=\K\otimes_\C(\F_{\al,1}\otimes\F^{R}\otimes e^{\bar{\Lambda}_l+\bar{\Lambda}_1+\gamma})\otimes e^{Q_{\bar{\mu}+\kappa}}.
\en
The highest weight vectors are given by $1\otimes  1\otimes 1\otimes e^{Q_{\bar{\mu}}}$ for $\hV(\Lambda_{0}+\mu,\mu)$, $1\otimes  1\otimes e^{\bar{\Lambda}_1}\otimes e^{Q_{\bar{\mu}}}$ for $\hV(\Lambda_{1}+\mu,\mu)$ and $1\otimes  1\otimes \left(\mmatrix{1\cr 1\cr}\right)\otimes e^{\bar{\Lambda}_l}\otimes e^{Q_{\bar{\mu}}}$ for $\hV(\Lambda_{l}+\mu,\mu)$,   
respectively. 

It is also easy to calculate the characters of these modules:
\be
&&ch_{\hV(\Lambda_{a}+\mu,\mu)}=\tr_{\hV(\Lambda_{a}+\mu,\mu)}q^{-{d}-\frac{c_W}{24}}=\sum_{\la\in {max}(\Lambda_a)\atop {mod\ Q_0+\C\delta\atop \gamma\in \cQ_0,\kappa\in \cQ}}ch_{\F_{\la,\gamma,\kappa}(\Lambda_a,\mu)},
\en
where $c_{W}=(l+\frac{1}{2})\left(1-\frac{2l(2l-1)}{rr^*}\right)$ is the central charge of the $W\! B_l$ algebra by Fateev and Lukyanov\cite{FaLu}, and 
the derivation operator $\hd$ is realized as  
\be
\hd&=&-\frac{1}{2}\sum_{j=1}^lh_jh^j-N^\al-N^\Psi+\frac{1}{2r^*}\sum_{j=1}^l(P_j+2)P^j-\frac{1}{2r}\sum_{j=1}^l((P+h)_j+2)(P+h)^j,\\
N^\al&=&\sum_{j=1}^l\sum_{m\in\Z_{>0}}\frac{m^2}{[m]}\frac{1-p^{*m}}{1-p^m} q^m\al_{j,-m}A^j_m,\quad
N^\Psi=\sum_{m{>0}}\frac{m(q^{\frac{1}{2}}+q^{-\frac{1}{2}})}{q^m+q^{-m}}\Psi_{-m}\Psi_m
\en
where $r, r^*\in \C^\times$, and  $A^j_m$ are the fundamental weight type elliptic bosons of the type $B_l$ given in Sec.5.1, $\Psi_m$ denotes $\Psi_m$ on $\F^{NS}$ and $\widehat{\Psi}_m$ on $\F^R$. We obtain:
\be
&&ch_{\hV(\Lambda_{a}+\mu,\mu)}=\sum_{\la\in {max}(\Lambda_a)\atop {mod\ Q_0+\C\delta\atop \gamma\in \cQ_0,\kappa\in \cQ}}ch_{\F_{\la,\gamma,\kappa}(\Lambda_a,\mu)},\\
&&ch_{\F_{\Lambda_0,\gamma,\kappa}(\Lambda_0,\mu)}=c^{\Lambda_0}_{\Lambda_0} 
q^{\frac{1}{2rr^*}|r(\bar{\mu}+\kappa+\bar{\rho})-r^*(\bar{\mu}+\kappa+\gamma+\bar{\rho})|^2},\\
&&ch_{\F_{\Lambda_1,\gamma,\kappa}(\Lambda_0,\mu)}=c^{\Lambda_0}_{\Lambda_1} 
q^{\frac{1}{2rr^*}|r(\bar{\mu}+\kappa+\bar{\rho})-r^*(\bar{\Lambda}_1+\bar{\mu}+\kappa+\gamma+\bar{\rho})|^2},\\
&&ch_{\F_{\Lambda_1,\gamma,\kappa}(\Lambda_1,\mu)}=c^{\Lambda_1}_{\Lambda_1} 
q^{\frac{1}{2rr^*}|r(\bar{\mu}+\kappa+\bar{\rho})-r^*(\bar{\Lambda}_1+\bar{\mu}+\kappa+\gamma+\bar{\rho})|^2},\\
&&ch_{\F_{\Lambda_0,\gamma,\kappa}(\Lambda_1,\mu)}=c^{\Lambda_1}_{\Lambda_0} 
q^{\frac{1}{2rr^*}|r(\bar{\mu}+\kappa+\bar{\rho})-r^*(\bar{\mu}+\kappa+\gamma+\bar{\rho})|^2},\\
&&ch_{\F_{\Lambda_l,\gamma,\kappa}(\Lambda_l,\mu)}=c^{\Lambda_l}_{\Lambda_l} 
q^{\frac{1}{2rr^*}|r(\bar{\mu}+\kappa+\bar{\rho})-r^*(\bar{\Lambda}_l+\bar{\mu}+\kappa+\gamma+\bar{\rho})|^2},\\
&&ch_{\F_{\Lambda_l-\al_l,\gamma,\kappa}(\Lambda_1,\mu)}=c^{\Lambda_1}_{\Lambda_l-\al_l} 
q^{\frac{1}{2rr^*}|r(\bar{\mu}+\kappa+\bar{\rho})-r^*(\bar{\Lambda}_l+\bar{\mu}+\kappa+\gamma+\bar{\rho})|^2},
\en
where 
\be
&&c^{\Lambda_0}_{\Lambda_0}=c^{\Lambda_1}_{\Lambda_1}=
\frac{q^{-\frac{1}{48}}}{2\eta(q)^l}\left((-q^{\frac{1}{2}};q)_\infty+(q^{\frac{1}{2}};q)_\infty\right),\\
&&c^{\Lambda_1}_{\Lambda_0}=c^{\Lambda_0}_{\Lambda_1}=\frac{q^{-\frac{1}{48}}}{2\eta(q)^l}
\left((-q^{\frac{1}{2}};q)_\infty-(q^{\frac{1}{2}};q)_\infty\right),\\
&&c^{\Lambda_l}_{\Lambda_l}=c^{\Lambda_l}_{\Lambda_l-\al_l}=\frac{q^{\frac{1}{24}}}{2\eta(q)^l}
(-q;q)_\infty. 
\en
$\sum_{\la\in {max}(\Lambda_a)\atop {mod\ Q_0+\C\delta}}ch_{\F_{\la,\gamma,\kappa}(\Lambda_a,\mu)}$ coincides with 
the character of the Verma modules of the $W\! B_l$-algebra with the highest weight 
$h=\frac{1}{2rr^*}|r(\bar{\mu}+\kappa+\bar{\rho})-r^*(\bar{\Lambda}_a+\bar{\mu}+\gamma+\kappa+\bar{\rho})|^2$ 
and the central charge $c_W$ with $r, r^*=r-1\in \C$ being generic. 

\begin{conj}
There exists a deformation of the $W\! B_l$-algebra such that 
\begin{itemize}
\item[i)] its generating functions commute with the level-$1$ elliptic currents  $e_j(z)$ and $f_j(z)$ of $U_{q,p}(B^{(1)}_l)$ modulo a total difference, i.e. $e_j(z)$ and $f_j(z)$ at $c=1$ are the  screening currents of the deformation of the  $W\! B_l$-algebra, 
\item[ii)]  for generic $r$ and $\mu\in \h^*$,  $\F_{\la,\xi,\kappa}(\Lambda+\mu,\mu)$ 
is an irreducible module of the deformation of the $W\! B_l$-algebra. 
\end{itemize}

\end{conj}

\noindent
{\it Remark.}\ All the algebras $W(\g)$ appearing in sec.\ref{ADE} and $W\! B_l$ in this subsection are the $W$-algebras  
associated with the coset $X_l^{(1)}\oplus X_l^{(1)}\supset (X_l^{(1)})_{{\rm diag}}$ 
with level $(r-g-1,1)$. In particular, the $W\! B_l$ is different from the one obtained from the quantum Hamiltonian reduction of the affine Lie algebra $B_l^{(1)}$. The $W$-algebras associated with such coset describe the critical behavior of the face type solvable lattice models introduced by Jimbo, Miwa and Okado\cite{JMO}.

\section{Elliptic Bosons of Various Types}
In this section we introduce elliptic bosons of the fundamental weight type $A^j_m$  
and the orthogonal basis type $\cE_m^{\pm j}$ 
 for $U_{q,p}(\gh)$, $\gh=A_l^{(1)}, B_l^{(1)}, C_l^{(1)}, D_l^{(1)}$. 
The level-1 bosons $A^j_m$ and $\cE_m^{\pm j}$ are used to realize the derivation operator $d$ and the 
 generating function of the deformed $W(A_l)$-algebra, respectively, in sec.\ref{level1Ex}. 

\subsection{Definition}\lb{fundamentalboson}
Let us set $\eta=-{tg}/2\ (t=(\mbox{long root})^2/2)$. 
$$
\begin{array}{c|cccc}
&A_l^{(1)}&B_l^{(1)}&C_l^{(1)}&D_l^{(1)}\\ \hline
g&l+1 &2l-1 &n+1 &2l-2\\
t&1&1&2&1\\
\eta&-\frac{l+1}{2}&-\frac{2l-1}{2}&-(l+1)&-(l-1)\\
\end{array}
$$
 Let $\al_{i,m}$ be the elliptic bosons of the simple root type as in Sec.2. 
We define the fundamental weight type elliptic bosons $A^j_m\ (1\leq j\leq l, m\in \Z_{\not=0})$ 
by
\bea
[\al_{i,m},A^j_n]=-\delta_{i,j}\delta_{m+n,0}\frac{[cm]}{m}\frac{1-p^m}{1-p^{*m}}q^{-cm} \qquad (1\leq i,j\leq l).\lb{alA}
\ena
Note that using the matrix $B(m)=([b_{i,j}m])_{1\leq i,j\leq l}$, we have\cite{FrRe}
\be
&&A^j_m=\sum_{k=1}^l(B(m)^{-1})_{kj}\al_{k,m}.
\en 
Solving \eqref{alA} we obtain the following. 

\noindent
For $A_l^{(1)}$, 
\be
&&A^j_m=
C_m\left(
[(2\eta+j)m]\sum_{k=1}^{j}[km]\al_{k,m}
+[jm]\sum_{k=j+1}^l[(2\eta+k)m]
\al_{k,m}\right)\qquad (1\leq j\leq l).
\en

\noindent
For $B_l^{(1)}$, 
\be
&&A^j_m=
C_m\left((q^{(\eta+j)m}+q^{-(\eta+j)m})
\sum_{k=1}^{j}[km]\al_{k,m}
+[jm]\sum_{k=j+1}^l(q^{(\eta+k)m}+q^{-(\eta+k)m})
\al_{k,m}\right)\\
&&\qquad\qquad\qquad\qquad\qquad (1\leq j\leq l).
\en

\noindent 
For $C_l^{(1)}$, 
\be
&&A^{j}_{m}=C_m \left((q^{(\eta+j)m}+q^{-(\eta+j)m})\sum_{k=1}^{j} [km]\alpha_{k,m}\right.\\
&&\qquad\quad
+ \left.[jm]\sum_{k=j+1}^{l-1}(q^{(\eta +k)m}+q^{-(\eta +k)m})\alpha_{k,m}+[jm]\alpha_{l,m} \right),
\quad(1\leq j\leq l-1),\\
&&A^{l}_{m}=C_m 
\left(\sum_{k=1}^{l-1} [km]\alpha_{k,m}+\frac{[m]}{[2m]}[lm]\alpha_{l,m} \right).
\en

\noindent
For $D_l^{(1)}$, 
\be
A^{j}_m&=&C_m\left(
(q^{(\eta+j)m}+q^{-(\eta+j)m})\sum_{k=1}^{j}[km]\al_{k,m}\right.\\
&&\qquad\qquad \left.
+[jm]\sum_{k=j+1}^{l-2}(q^{(\eta+k)m}+q^{-(\eta+k)m})\al_{k,m}+[jm](\al_{l-1,m}+a_{l,m})\right)\quad (1\leq j\leq l-2),\\
A^{l-1}_m&=&C_m\left(\sum_{k=1}^{l-2}[km]\al_{k,m}
+\frac{[m]}{[2m]}([lm]\al_{l-1,m}+[(l-2)m]a_{l,m})\right),\\
A^{l}_m&=&C_m\left(\sum_{k=1}^{l-2}[km]\al_{k,m}
+\frac{[m]}{[2m]}([(l-2)m]\al_{l-1,m}+[lm]a_{l,m})\right).
\en

Here
\be
C_m&=&\frac{1}{[m]^2[2\eta m]}\qquad\mbox{for} A_l^{(1)}\\
&=&\frac{[\eta m]}{[m]^2[2\eta m]}\qquad\mbox{for} B_l^{(1)}, C_l^{(1)}, D_l^{(1)}.
\en

We then devide $A^j_m$ into two terms and define the elliptic bosons $\cE^{\pm j}_m$ of the orthogonal basis type as follows. 

\noindent
For $A_l^{(1)}$,
\bea
&&A^j_m=\cE^{+j}_m+\cE^{-j}_m,\\
&&\cE^{\pm j}_m=\pm q^{\pm jm}\frac{C_m}{q-q^{-1}}\left(q^{\pm 2\eta m}\sum_{k=1}^{j-1}[km]\al_{k,m}
+\sum_{k=j}^l [ (2\eta+k ) m]\al_{k,m} \right) \lb{ONboson}
\ena
for $1\leq j\leq l$.
It is convenient to define $\cE^{\pm (l+1)}_m$ by
\bea
&&\cE^{\pm (l+1)}_m=\mp\frac{C_m}{q-q^{-1}}\sum_{k=1}^l[km]\al_{k,m}. 
\ena

\noindent
For $B_l^{(1)}$, 
\bea
&&A^j_m=\cE^{+j}_m+\cE^{-j}_m,\\
&&\cE^{\pm j}_m=q^{\pm jm}C_m\left(q^{\pm \eta  m}\sum_{k=1}^{j-1}[km]\al_{k,m}
\pm\sum_{k=j}^l [ ( \eta+k ) m]_+\al_{k,m} \right) 
\ena
for $1\leq j\leq l$. Here we set
\be
[m]_+=\frac{q^m+q^{-m}}{q-q^{-1}}.
\en
We also define 
\bea
\cE^0_m=\frac{[\frac{m}{2}]}{[m]}(\cE^{+l}_m+\cE^{-l}_m).
\ena

\noindent 
For $C_l^{(1)}$, 
\bea
&&A^{j}_{m}=\cE_m^{+j}+\cE_m^{-j}\\
&& \cE^{\pm j}_{m}=q^{\pm j m}C_{m}\left( q^{\pm\eta m}\sum_{k=1}^{j-1} [km]\alpha_{k,m}\pm \sum_{k=j}^{l-1}
[(\eta+k)m]_+\alpha_{k,m} \pm \frac{\alpha_{l,m}}{q-q^{-1}}
 \right)\qquad(1\leq j\leq l-1),\nn\\
&&\\
&&A_m^{l}=\frac{1}{q^{m}+q^{-m}}(\cE_m^{+l}+\cE_m^{-l}),\\
&& \cE^{\pm l}_{m}=q^{\pm lm}C_m\left( q^{\pm \eta m}\sum_{k=1}^{l-1} [km]\alpha_{k,m}\pm\frac{\alpha_{l,m}}{q-q^{-1}} \right). 
\ena

\noindent
For $D_l^{(1)}$, 
\bea
&&A^j_m=\cE^{+j}_m+\cE^{-j}_m,\\
&&\cE^{\pm j}_m=q^{\pm jm}C_m\left(q^{\pm \eta m}\sum_{k=1}^{j-1}[km]\al_{k,m}
\pm \sum_{k=j}^{l-2} [(\eta+k)m]_+\al_{k,m} \pm \frac{1}{q-q^{-1}}(\al_{l-1,m}+\al_{l,m})\right)\nn\\
&&\qquad\qquad\qquad\qquad (1\leq j\leq l-2),\\
&&\cE^{\pm (l-1)}_m=C_m\left(\sum_{k=1}^{l-2}[km]\al_{k,m}
 \pm \frac{q^{\mp\eta m}}{q-q^{-1}}(\al_{l-1,m}+\al_{l,m})\right),\\
&&\cE^{\pm l}_m=q^{\pm m}C_m\left(\sum_{k=1}^{l-2}[km]\al_{k,m}
 \mp \frac{1}{q-q^{-1}}(q^{\pm\eta m}\al_{l-1,m}-q^{\mp\eta m}\al_{l,m})\right).
\ena

\begin{prop}\lb{alphaAA}
\noindent 
\bea
&&\al_{j,m}=\pm[m]^2(q-q^{-1})(\cE^{\pm j}_m-q^{\mp m}\cE^{\pm(j+1)}_m), 
\ena
$1\leq j\leq l$ for $A_l^{(1)}$, $1\leq j\leq l-1$ for $B_l^{(1)}$, $C_l^{(1)}$, $D_l^{(1)}$,  and 
\bea
\al_{l,m}&=&[m](q^{m/2}-q^{-m/2})(q^{-m/2}\cE^{+l}_m-q^{m/2} \cE^{-l}_m) \quad \mbox{for\ } B_l^{(1)}, \\
&=&[m]^2(q-q^{-1})\left(q^{m}\cE^{+ l}_m-q^{-m}\cE^{-l}_m\right)\quad \mbox{for\ } C_l^{(1)}, \\
&=&\pm[m]^2(q-q^{-1})(\cE^{\pm (l-1)}_m-q^{\pm m}\cE^{\mp l}_m)\quad \mbox{for\ } D_l^{(1)}.
\ena
\end{prop}

\begin{prop}\lb{cEA}
The following relations hold.
\bea
&&\cE^{\pm 1}_m=\pm\frac{q^{\pm m}}{q^m-q^{-m}}A^{1}_m,\qquad \cE^{\pm j}_m=\pm\frac{1}{q^m-q^{-m}}\left(q^{\pm m}A^{j}_m-A^{j-1}_m\right), 
\ena
where  $2\leq j\leq l$ for $A_l^{(1)}$, 
$2\leq j\leq l$ for $B_l^{(1)}$, $2\leq j\leq l-1$ for  $C_l^{(1)}$ 
and $2\leq j\leq l-2$ for $D_l^{(1)}$.  In addition, we have 
\bea
&&\cE^{\pm (l+1)}_m=\mp\frac{1}{q^m-q^{-m}}A^{l}_m,\qquad \sum_{j=1}^{l+1}q^{\pm(j-1)m}\cE^{\pm j}_m=0 
\qquad \mbox{for}\ A_l^{(1)}, \\
&&\cE^{\pm l}_m=\pm\frac{1}{q^m-q^{-m}}\left((q^m+q^{-m})q^{\pm m}A^l_m-A^{l-1}_m\right)
\qquad \mbox{for}\ C_l^{(1)},
\ena
and
\bea
&&\cE^{\pm(l-1)}_m=\pm\frac{1}{q^m-q^{-m}}\left(q^{\pm m}A^{l-1}_m+q^{\pm m}A^l_m-A^{l-2}_m\right),\\
&&\cE^{\pm l}_m=\pm\frac{1}{q^m-q^{-m}}\left(q^{\pm 2m}A^l_m-A^{l-1}_m\right)\qquad\qquad \mbox{for}\ D_l^{(1)}.
\ena
\end{prop}
\noindent
{\it Remark.}\ 
The level-1 case i.e. $c=1$, the $A^{(1)}_l$ type relation was given in \cite{FeFr, AKOS} and  the $D^{(1)}_l$ type was essentially given in \cite{FrRe}, where parameters $q$ and $t$ should be identified with our $p^{*\frac{1}{2}}=p^{\frac{1}{2}}q^{-1}$ and $p^{\frac{1}{2}}$, respectively. However the $B^{(1)}_l$ and $C^{(1)}_l$ cases are  different from those given in \cite{FrRe}. At least the formulas for $B_l^{(1)}$ and $C^{(1)}_l$ seem to be reversed. Our definitions and relations are valid for arbitrary level $c$. 

Although the expressions of $\cE^{\pm j}_m$ are complicated depending on the types of the affine Lie algebras, 
their commutation relations are  rather  universal:
\begin{thm}\lb{cEcE}
For $1\leq j, k \leq l$, the following commutation relations hold. 
For $A_l^{(1)}$, 
\bea
&&[\cE^{\pm j}_m,\cE^{\pm j}_n]=[\cE^{\pm j}_m,\cE^{\mp j}_n]
=\delta_{m+n,0}\frac{[cm][(2\eta+1) m] }{m(q-q^{-1})^2[m]^3 [2\eta m]}\frac{1-p^m}{1-p^{*m}}q^{-cm}, \\
&&[\cE^{\pm j}_m,\cE^{\pm k}_n]
=\delta_{m+n,0}q^{\mp({\rm sgn}(k-j)2\eta +k-j)m}
\frac{[cm]}{m(q-q^{-1})[m]^2[2\eta m]}\frac{1-p^m}{1-p^{*m}}q^{-cm},\\
&&
[\cE^{\pm j}_m,\cE^{\mp k}_n]=- \delta_{m+n,0}q^{\pm(2\eta+j+k)m}\frac{[cm]}{m(q-q^{-1})[m]^2[2\eta m]}\frac{1-p^m}{1-p^{*m}}q^{-cm}{.}
\ena
For $B_l^{(1)}, C_l^{(1)}, D_l^{(1)}$, 
\bea
&&[\cE^{\pm j}_m,\cE^{\pm j}_n]
=\delta_{m+n,0}\frac{[cm][\eta m] [2(\eta+1)m]}{m(q-q^{-1})^2[m]^3[2\eta m] [(\eta+1)m]}\frac{1-p^m}{1-p^{*m}}q^{-cm}, 
\\
&&[\cE^{\pm j}_m,\cE^{\mp j}_n]
=\mp\delta_{m+n,0}\frac{q^{\pm jm}[cm][\eta m]}{m[m]^3(q-q^{-1})^2[2\eta m]}
\frac{1-p^m}{1-p^{*m}}q^{-cm} 
\left(q^{\pm(\eta +j)m}[m]\pm q^{\mp(j-1)m}[\eta m]_{+}\right),\nn\\
&&\\
&&[\cE^{\pm j}_m,\cE^{\pm k}_n]
=\mp{\rm sgn}(k-j)\delta_{m+n,0}q^{\mp({\rm sgn}(k-j) \eta +k-j)m}
\frac{[cm][\eta m]}{m(q-q^{-1})[m]^2[2 \eta m]}\frac{1-p^m}{1-p^{*m}}q^{-cm},\nn\\
&&\\
&&
[\cE^{\pm j}_m,\cE^{\mp k}_n]=\mp \delta_{m+n,0}q^{\pm(\eta+j+k)m}\frac{[cm][\eta m]}{m(q-q^{-1})[m]^2[2 \eta m]}\frac{1-p^m}{1-p^{*m}}q^{-cm}{.}
\ena
Here
\be
{\rm sgn}(l-j)=\left\{\mmatrix{+&(l>j),\cr
                              -&(l<j). \cr }\right.
\en
\end{thm}
\noindent
{\it Proof.} 
Straightforward calculation using  Proposition \ref{cEA} and \eqref{alA}.  
\qed
\begin{prop}
For $1\leq i\leq l$,   the following commutation relations hold. 
\bea
&&[\al_{i,m}, \cE^{\pm j}_n]=\pm \delta_{m+n,0}\frac{[cm]}{m(q^m-q^{-m})}\frac{1-p^m}{1-p^{*m}}q^{-cm}(q^{\mp m}\delta_{i,j}-\delta_{i,j-1})
\ena
where  $1\leq j\leq l$ for $A^{(1)}_l, B_l^{(1)}$,  
$1\leq j\leq l-1$ for $C_l^{(1)}$, $1\leq j\leq l-2$ for $D_l^{(1)}$. 
In addition, 
\bea
&&[\al_{i,m}, \cE^{\pm l}_n]=\pm \delta_{m+n,0}\frac{[cm]}{m(q^m-q^{-m})}\frac{1-p^m}{1-p^{*m}}q^{-cm}(q^{\mp m}(q^m+q^{-m})\delta_{i,l}-\delta_{i,l-1})\ \mbox{for}\ C_l^{(1)}{,}\nn\\
&&
\ena  
 and  
\bea
&&[\al_{i,m}, \cE^{\pm (l-1)}_n]=\pm \delta_{m+n,0}\frac{[cm]}{m(q^m-q^{-m})}\frac{1-p^m}{1-p^{*m}}q^{-cm}(q^{\mp m}\delta_{i,l-1}+q^{\mp m}\delta_{i,l}-\delta_{i,l-2}),\nn\\
&&\\
&&[\al_{i,m}, \cE^{\pm l}_n]=\pm \delta_{m+n,0}\frac{[cm]}{m(q^m-q^{-m})}\frac{1-p^m}{1-p^{*m}}q^{-cm}(q^{\mp 2m}\delta_{i,l}-\delta_{i,l-1})\qquad \mbox{for}\ D_l^{(1)}.
\ena 
\end{prop}

From \eqref{bosonve} and \eqref{bosonvf} we also obtain the following relations. 
\begin{prop}\lb{cEef}
For $1\leq j\leq l$, 
\bea
&&[\cE^{\pm i}_m, e_j(z)]=\pm\frac{q^{-cm}z^m}{m(q^m-q^{-m})}\frac{1-p^m}{1-p^{*m}}e_j(z)(q^{\pm m}\delta_{i,j}-\delta_{i-1,j}),\\
&&[\cE^{\pm i}_m, f_j(z)]=\mp\frac{z^m}{m(q^m-q^{-m})}f_j(z)(q^{\pm m}\delta_{i,j}-\delta_{i-1,j})
\ena
where  $1\leq i\leq l$ for $A^{(1)}_l, B_l^{(1)}$, $1\leq i\leq l-1$ for $C_l^{(1)}$, $1\leq i\leq l-2$ for $D_l^{(1)}$.
In addition,  
\bea
&&[\cE^{\pm l}_m, e_j(z)]=\pm\frac{q^{-cm}z^m}{m(q^m-q^{-m})}\frac{1-p^m}{1-p^{*m}}e_j(z)(q^{\pm m}(q^m+q^{-m})\delta_{l,j}-\delta_{l-1,j}),\\
&&[\cE^{\pm l}_m, f_j(z)]=\mp\frac{z^m}{m(q^m-q^{-m})}f_j(z)(q^{\pm m}(q^m+q^{-m})\delta_{l,j}-\delta_{l-1,j})\quad \mbox{for}\ C_l^{(1)}{,}
\ena
and 
\bea
&&[\cE^{\pm (l-1)}_m, e_j(z)]=\pm\frac{q^{-cm}z^m}{m(q^m-q^{-m})}\frac{1-p^m}{1-p^{*m}}e_j(z)(q^{\pm m}\delta_{l-1,j}+q^{\pm m}\delta_{l,j}-\delta_{l-2,j}),\\
&&[\cE^{\pm (l-1)}_m, f_j(z)]=\mp\frac{z^m}{m(q^m-q^{-m})}f_j(z)
(q^{\pm m}\delta_{l-1,j}+q^{\pm m}\delta_{l,j}-\delta_{l-2,j}),\\
&&[\cE^{\pm l}_m, e_j(z)]=\pm\frac{q^{-cm}z^m}{m(q^m-q^{-m})}\frac{1-p^m}{1-p^{*m}}e_j(z)(q^{\pm 2m}\delta_{l,j}-\delta_{l-1,j}),\\
&&[\cE^{\pm l}_m, f_j(z)]=\mp\frac{z^m}{m(q^m-q^{-m})}f_j(z)(q^{\pm 2m}\delta_{l,j}-\delta_{l-1,j})
\qquad \mbox{for}\ D_l^{(1)}.
\ena
\end{prop}

\subsection{The Elliptic Currents $k_{\pm j}(z)$}
Let us set 
\bea
&&\psi_j(z)=:\exp\left\{(q-q^{-1})\sum_{m\not=0}\frac{\al_{j,m}}{1-p^m}p^mz^{-m}\right\}:.
\ena
Then the elliptic currents $\psi_j^\pm(z)$ in Definition \ref{defUqp} can be written as 
\bea
&&\psi^+_j(q^{-\frac{c}{2}}z)=K^+_j\psi_j(z),\qquad \psi^-_j(q^{-\frac{c}{2}}z)=K^-_j\psi_j(pq^{-c}z).
\ena
Let us introduce the new currents $k_{\pm j}(z) \,  (1\leq j\leq l)$ associated with $\cE^{\pm j}_m$ by
\bea
k_{\pm j}(z) &=& :\exp\left\{\sum_{m\not=0} \frac{[m]^2(q-q^{-1})^2}{1-p^m}p^m\cE^{\pm j}_mz^{-m} \right\}: 
\ena 
and in addition we define $k_0(z)$ for $B_l^{(1)}$ by 
\bea
k_0(z) &=&:k_{-l}(q^{-1/2}z)\psi_l(q^{-1/2}z):=:k_{+l}(q^{1/2}z)\psi_l(q^{1/2}z)^{-1}:.
\ena
Then from Proposition \ref{alphaAA} we have the following decompositions.  
\begin{prop}\lb{psikk}
\bea
\psi_j(z)&=&:k_{+j}(z)k_{+(j+1)}(qz)^{-1}:=:k_{-j}(z)^{-1}k_{-(j+1)}(q^{-1}z):
\ena
where $1\leq j\leq l-1$ for  $A_l^{(1)}$, $1\leq j\leq l-1$ for $B_l^{(1)}$, $C_l^{(1)}$ and $D_l^{(1)}$. 
In addition, 
\bea
\psi_l(z)&=&:k_{+l}(z)k_0(q^{-1/2}z)^{-1}:=:k_{-l}(z)^{-1}k_0(q^{1/2}z): \qquad \mbox{for}\ B_l^{(1)},\\
&=&:k_{+ l}(q^{-1}z) k_{- l}(q z)^{- 1}:\qquad\qquad\qquad\qquad\qquad\qquad \mbox{for}\ C_l^{(1)},\\  
&=&:k_{+(l-1)}(z)k_{-l}(q^{-1}z)^{-1}:=:k_{-(l-1)}(z)^{-1}k_{+l}(qz):\quad \mbox{for}\ D_l^{(1)}.
\ena
\end{prop}

Now let us introduce the functions $\tilde\rho^+(z)$, which appear associated with 
the elliptic dynamical $R$-matrices\cite{Konno06}:
\bea
\tilde{\rho}^+(z)&=&\frac{\{q^2 z\}\{\xi^2q^{-2} z\}}{\{\xi^2 z\}\{z\}}
\frac{\{p\xi^2/z\}\{p/z\}}{\{p\xi^2q^{-2}/z\}\{pq^2/z\}}\qquad \mbox{for\ }A^{(1)}_l,\\
&=&\frac{\{\xi z\}^{2}\{\xi^2q^{-2} z\}\{q^2 z\}}{\{\xi^2 z\}\{z\}\{\xi q^2 z\}\{\xi q^{-2} z\}}
\frac{\{p\xi^2/z\}\{p/z\}\{p\xi q^2/z\}\{p\xi q^{-2}/z\}}{\{p\xi/z\}^{2}\{p\xi^2q^{-2}/z\}\{pq^2/z\}}\quad \mbox{for\ }B^{(1)}_l, C^{(1)}_l, D^{(1)}_l, \nn\\
&&
\ena
where $\xi=q^{-2\eta}$, $\{z\}=(z;p,\xi^2)_\infty$.  
The following Theorem indicates a deep relationship between  $k_{\pm j}(z)$'s and elliptic dynamical $R$-matrices. 
\begin{thm}\lb{kk}
\be
&&k_{\pm j}(z_1)k_{\pm j}(z_2)=\frac{\tilde{\rho}^{+*}(z)}{\tilde{\rho}^+(z)}k_{\pm j}(z_2)k_{\pm j}(z_1),
\qquad  (1\leq j\leq l),\\
&&k_{+j}(q^jz_1)k_{+k}(q^kz_2)=\frac{\tilde{\rho}^{+*}(z)}{\tilde{\rho}^+(z)}\frac{\Theta_{p^*}(q^{-2}z)\Theta_{p}(z)}{\Theta_{p^*}(z)\Theta_{p}(q^{-2}z)} k_{+k}(q^kz_2)k_{+j}(q^jz_1) \qquad \qquad (1\leq j<k\leq l),\\
&&k_{-j}(q^{-j}z_1)k_{-k}(q^{-k}z_2)=\frac{\tilde{\rho}^{+*}(z)}{\tilde{\rho}^{+}(z)}\frac{\Theta_{p^*}(q^{-2}z)\Theta_{p}(z)}{\Theta_{p^*}(z)\Theta_{p}(q^{-2}z)} k_{-k}(q^{-k}z_2)k_{-j}(q^{-j}z_1) \qquad (1\leq k<j\leq l),\\
&&k_{+j}(q^{j}z_1)k_{-k}(q^{-k}\xi z_2)=\frac{\tilde{\rho}^{+*}(z)}{\tilde{\rho}^{+}(z)}\frac{\Theta_{p^*}(q^{-2}z)\Theta_{p}(z)}{\Theta_{p^*}(z)\Theta_{p}(q^{-2}z)} k_{-k}(q^{-k}\xi z_2)k_{+j}(q^{j}z_1) \qquad (j\not=k),\\
&&k_{+j}(q^jz_1)k_{-j}(q^{-j}\xi z_2)=\frac{\tilde{\rho}^{+*}(u)}{\tilde{\rho}^{+}(u)}
\frac{\Theta_{p^*}(q^{2j-2}\xi^{-1}z)\Theta_{p}(q^{2j}\xi^{-1}z)}{\Theta_{p^*}(q^{2j}\xi^{-1}z)\Theta_{p}(q^{2j-2}\xi^{-1}z)} \frac{\Theta_{p^*}(q^{-2}z)\Theta_{p}(z)}{\Theta_{p^*}(z)\Theta_{p}(q^{-2}z)} k_{-j}(q^{-j}\xi z_2)k_{+j}(q^jz_1),
\en
where $z=z_1/z_2$ and   $\tilde{\rho}^{+*}(z)=\tilde{\rho}^+(z)|_{p\mapsto p^*}$. 
In addition, for $B_l^{(1)}$ we have 
\be
&&k_0(z_1)k_0(z_2)=\frac{\tilde{\rho}^{+*}(u)}{\tilde{\rho}^{+}(u)}
\frac{\Theta_{p^*}(q^{-2}z)\Theta_{p}(q^{2}z)\Theta_{p^*}(qz)\Theta_{p}(q^{-1}z)}
{\Theta_{p^*}(q^{2}z)\Theta_{p}(q^{-2}z)\Theta_{p^*}(q^{-1}z)\Theta_{p}(qz)}
k_0(z_2)k_0(z_1),\\
&&k_{+j}(q^jz_1)k_0(q^{l-1/2}z_2)=\frac{\tilde{\rho}^{+*}(u)}{\tilde{\rho}^{+}(u)}
\frac{\Theta_{p^*}(q^{-2}z)\Theta_{p}(z)}
{\Theta_{p^*}(z)\Theta_{p}(q^{-2}z)}
k_0(q^{l-1/2}z_2)k_{+j}(q^jz_1)\quad (1\leq j\leq l),\\
&&k_{-j}(\xi q^{-j}z_1)k_0(q^{l-1/2}z_2)=\frac{\tilde{\rho}^{+*}(u)}{\tilde{\rho}^{+}(u)}
\frac{\Theta_{p^*}(z)\Theta_{p}(q^{2}z)}
{\Theta_{p^*}(q^{2}z)\Theta_{p}(z)}
k_0(q^{l-1/2}z_2)k_{-j}(\xi q^{-j}z_1)\quad (1\leq j\leq l).
\en
\end{thm}
\noindent
{\it Proof.} 
Straightforward calculation using Theorem \ref{cEcE}.  
\qed

In addition from Proposition \ref{cEef}, we obtain:
\begin{prop}\lb{kef}
\be
&&k_{\pm j}(z_1)e_j(z_2)=\frac{\Theta_{p^*}(q^{-c}z)}{\Theta_{p^*}(q^{-c\mp 2}z)}
e_j(z_2)k_{\pm j}(z_1)\qquad (1\leq j\leq l),\\
&&k_{\pm j}(z_1)e_{j-1}(z_2)=\frac{\Theta_{p^*}(q^{-c\mp 1}z)}{\Theta_{p^*}(q^{-c\pm 1}z)}
e_{j-1}(z_2)k_{\pm j}(z_1)\qquad (2\leq j\leq l),\\
&&k_{\pm j}(z_1)e_k(z_2)=e_k(z_2)k_{\pm j}(z_1)\qquad(k\not=j,j-1), 
\en
\be
&&k_{\pm j}(z_1)f_j(z_2)=\frac{\Theta_{p}(q^{\mp 2}z)}{\Theta_{p}(z)}
f_j(z_2)k_{\pm j}(z_1)\qquad (1\leq j\leq l),\\
&&k_{\pm j}(z_1)f_{j-1}(z_2)=\frac{\Theta_{p}(q^{\pm 1}z)}{\Theta_{p}(q^{\mp 1}z)}
f_{j-1}(z_2)k_{\pm j}(z_1)\qquad (2\leq j\leq l),\\
&&k_{\pm j}(z_1)f_k(z_2)=f_k(z_2)k_{\pm j}(z_1)\qquad(k\not=j,j-1)
\en
for $A^{(1)}_l, B_l^{(1)}$ with $1\leq i\leq l$, $C_l^{(1)}$ with $1\leq i\leq l-1$, $D_l^{(1)}$ with $1\leq i\leq l-2$.
In addition, we have
\be
&&k_0(q^{l-1/2}z_1)e_l(z_2)=\frac{\Theta_{p^*}(q^{-c+l}z)\Theta_{p^*}(q^{-c+l-1}z)}
{\Theta_{p^*}(q^{-c+l-2}z)\Theta_{p^*}(q^{-c+l+1}z)}e_l(z_2)k_0(q^{l-1/2}z_1),\\
&&k_0(q^{l-1/2}z_1)e_j(z_2)=e_j(z_2)k_0(q^{l-1/2}z_1)\qquad (1\leq j\leq l-1),\\
&&k_0(q^{l-1/2}z_1)f_l(z_2)=\frac{\Theta_{p}(q^{l-2}z)\Theta_{p}(q^{l+1}z)}
{\Theta_{p}(q^{l}z)\Theta_{p}(q^{l-1}z)}f_l(z_2)k_0(q^{l-1/2}z_1),\\
&&k_0(q^{l-1/2}z_1)f_j(z_2)=f_j(z_2)k_0(q^{l-1/2}z_1)\quad (1\leq j\leq l-1)\qquad\qquad \mbox{for}\  B_l^{(1)}, \\
&&k_{\pm l}(z_1)e_l(z_2)=\frac{\Theta_{p^*}(q^{-c\pm 1}z)}{\Theta_{p^*}(q^{-c\mp 3}z)}
e_l(z_2)k_{\pm l}(z_1),\\
&&k_{\pm l}(z_1)e_{l-1}(z_2)=\frac{\Theta_{p^*}(q^{-c\mp 1}z)}{\Theta_{p^*}(q^{-c\pm 1}z)}
e_{l-1}(z_2)k_{\pm l}(z_1),\\
&&k_{\pm l}(z_1)e_j(z_2)=e_j(z_2)k_{\pm l}(z_1)\qquad(j\not=l,l-1), \\
&&k_{\pm l}(z_1)f_l(z_2)=\frac{\Theta_{p}(q^{\mp 3}z)}{\Theta_{p}(q^{\pm 1}z)}
f_l(z_2)k_{\pm l}(z_1),\\
&&k_{\pm l}(z_1)f_{l-1}(z_2)=\frac{\Theta_{p}(q^{\pm 1}z)}{\Theta_{p}(q^{\mp 1}z)}
f_{l-1}(z_2)k_{\pm l}(z_1),\\
&&k_{\pm l}(z_1)f_j(z_2)=f_j(z_2)k_{\pm l}(z_1)\quad(j\not=l,l-1)\qquad\qquad \mbox{for}\  C_l^{(1)}, \\
&&k_{\pm (l-1)}(z_1)e_j(z_2)=\frac{\Theta_{p^*}(q^{-c}z)}{\Theta_{p^*}(q^{-c\mp 2}z)}
e_j(z_2)k_{\pm (l-1)}(z_1)\qquad (j=l,l-1),\\
&&k_{\pm (l-1)}(z_1)e_{l-2}(z_2)=\frac{\Theta_{p^*}(q^{-c\mp 1}z)}{\Theta_{p^*}(q^{-c\pm 1}z)}
e_{l-2}(z_2)k_{\pm (l-1)}(z_1),\\
&&k_{\pm (l-1)}(z_1)e_j(z_2)=e_j(z_2)k_{\pm (l-1)}(z_1)\qquad(j\not=l,l-1,l-2), \\
&&k_{\pm l}(z_1)e_l(z_2)=\frac{\Theta_{p^*}(q^{-c\mp 1}z)}{\Theta_{p^*}(q^{-c\mp 3}z)}
e_l(z_2)k_{\pm l}(z_1),\\
&&k_{\pm l}(z_1)e_{l-1}(z_2)=\frac{\Theta_{p^*}(q^{-c\mp 1}z)}{\Theta_{p^*}(q^{-c\pm 1}z)}
e_{l-1}(z_2)k_{\pm l}(z_1),\\
&&k_{\pm l}(z_1)e_j(z_2)=e_j(z_2)k_{\pm l}(z_1)\qquad(j\not=l,l-1), 
\\
&&k_{\pm (l-1)}(z_1)f_j(z_2)=\frac{\Theta_{p}(q^{\mp 2}z)}{\Theta_{p}(z)}
f_j(z_2)k_{\pm (l-1)}(z_1)\qquad (j=l,l-1),\\
&&k_{\pm (l-1)}(z_1)f_{l-2}(z_2)=\frac{\Theta_{p}(q^{\pm 1}z)}{\Theta_{p}(q^{\mp 1}z)}
f_{l-2}(z_2)k_{\pm (l-1)}(z_1),\\
&&k_{\pm (l-1)}(z_1)f_j(z_2)=f_j(z_2)k_{\pm (l-1)}(z_1)\qquad(j\not=l,l-1,l-2), 
\\
&&k_{\pm l}(z_1)f_l(z_2)=\frac{\Theta_{p}(q^{\mp 3}z)}{\Theta_{p}(q^{\mp 1}z)}
f_l(z_2)k_{\pm l}(z_1),\\
&&k_{\pm l}(z_1)f_{l-1}(z_2)=\frac{\Theta_{p}(q^{\pm 1}z)}{\Theta_{p}(q^{\mp 1}z)}
f_{l-1}(z_2)k_{\pm l}(z_1),\\
&&k_{\pm l}(z_1)f_j(z_2)=f_j(z_2)k_{\pm l}(z_1)\quad(j\not=l,l-1)\qquad\qquad\mbox{for}\  D_l^{(1)}.
\en
\end{prop}

The elliptic bosons $\cE^{\pm j}_m$ and their elliptic currents $k_{\pm j}(z)$ are useful to realize the $L$-operators 
and the vertex operators for $U_{q,p}(\gh)$ as well as deformation of the $W$-algebras. We will discuss this subject in  separate papers.

\section*{Acknowledgements}
H.K is supported by the Grant-in -Aid for Scientific Research (C) 22540022 JSPS, Japan.
R.M.F is grateful to the Egyptian government for a scholarship.

\begin{appendix}

\section{The Drinfeld Realization of $U_q(\gh)$}\lb{Uqgh}
Let $\gh$ be an untwisted affine Lie algebra. 

\begin{dfn}
The quantum affine algebra $U_{q}(\gh)$ in the Drinfeld realization is a unital $\C$-algebra 
generated by $q^h \ (h\in \h) $, $a^\vee_{i,n},\ x^\pm_{i, m}\ (i\in I,\ n\in \Z_{\not=0}, m\in \Z)$ $\bd$ and the central element $c$.   
We set 
\bea
&&x^\pm_i(z)=\sum_{m\in \Z}x^\pm_{i,m}z^{-m},\\
&&\psi_i(z)=q_i^{h_i}\exp\left((q_i-q_i^{-1})\sum_{n>0} a^\vee_{i,n}z^{-n}\right),\\
&&\varphi_i(z)=q_i^{-h_i}\exp\left(-(q_i-q_i^{-1})\sum_{n>0} a^\vee_{i,-n}z^{n}\right). 
\ena
The defining relations are as follows. 
\bea
&& [q_i^{\pm h_i},\bd]=0,\quad [\bd,a_{i,n}]=n a_{i,n},\quad 
[\bd,x^{\pm}_{i,n}]=n x^{\pm}_{i,n}, \lb{dxUq}\\
&&[q_i^{\pm h_i},a_{j,n}]=0,\qquad q_i^{h_i}x_j^\pm(z)=q_i^{\pm a_{ij}} x_j^{\pm}(z)q_i^{h_i},\lb{hxUq}\\
&&
[a^\vee_{i,n},a^\vee_{j,m}]=\frac{[a_{ij}n]_i[c n)_j}{n}q^{-c|n|}\delta_{n+m,0},\lb{aaUq}\\
&&
[a^\vee_{i,n},x_j^+(z)]=\frac{[a_{ij}n]_i}{n}q^{-c|n|}z^n x_j^+(z),\\
&&
[a^\vee_{i,n},x_j^-(z)]=-\frac{[a_{ij}n]_i}{n} z^n x_j^-(z),
\\
&&(z-q^{\pm b_{ij}}w)
x_i^\pm(z)x_j^\pm(w)= (q^{\pm b_{ij}}z-w) x_j^\pm(w)x_i^\pm(z),
\\
&&[x_i^+(z),x_j^-(w)]=\frac{\delta_{i,j}}{q_i-q_i^{-1}}
\left(\delta\bigl(q^{-k}\frac{z}{w}\bigr)\psi_i(q^{k/2}w)
-\delta\bigl(q^{k}\frac{z}{w}\bigr)\varphi_i(q^{-k/2}w)
\right),
\\
&&\sum_{\sigma\in S_a}\sum_{s=0}^a (-)^s
\left[\mmatrix{a\cr
s\cr}\right]_{i}
x^{\pm}_{i}(z_{\sigma(1)})\cdots x^{\pm}_{i}(z_{\sigma(s)})
x^{\pm}_{j}(w) x^{\pm}_{i}(z_{\sigma(s+l)})\cdots x^{\pm}_{i}(z_{\sigma(a)})=0,\nn\\
&&\qquad\quad (i\not=j,\ a=1-a_{ij}).\lb{serrefUq}
\ena
\end{dfn}

For $k\in \C$, we define the category ${C}_k$ of the level-$k$ $U_{q}(\gh)$-modules in the same way as  ${\cC}_k$ of 
$U_{q,p}(\gh)$ in Sec.2. 
Let $a_{i,n}=[d_i]a^\vee_{i,n} \ (i\in I, n\in \Z_{\not=0})$ be the simple root type level-$k$ Drinfeld bosons. 
They satisfy 
\be
&&[a_{i,n},a_{j,m}]=\frac{[b_{ij}n][kn]}{n}q^{-k|n|}\delta_{n+m,0}. 
\en
For $(V,\bar{\pi})\in {C}_k$,  we define the $Z$-operators associated with the level-$k$ $U_q(\gh)$-module 
$V$ by 
\be
&&{Z}^\pm_i(z;V)= \exp\left(\mp\sum_{n\geq 1}\frac{\bar{\pi}(a_{i,-n})}{[kn]}q^{\frac{1\mp 1}{2}kn}z^n\right)
\bar{\pi}(x_{i}^\pm(z))\exp\left(\pm\sum_{n\geq 1}\frac{\bar{\pi}(a_{i,n})}{[kn]}q^{\frac{1\mp 1}{2}kn}z^{-n}\right). \label{UqZ} 
\en
The coefficients $Z^\pm_{i,n}(V)$ of  ${Z}^\pm_i(z;V)=\sum_{n\in \Z}Z^\pm_{i,n}(V)z^{-n}$ in $z$ are well defined elements in $\End_\C V$. 
\begin{thm}
The $Z$-operators ${Z}^\pm_i(z;V)$ satisfy the same relations  in Theorem \ref{zz} 
except for \eqref{hZPZp},\eqref{hZPZm} with replacement $\cZ^\pm_j(z;\hV)$, $\al_{j,m}$, $d$ and $K^\pm_j$ by ${Z}^\pm_i(z;V), a_{j,m}$, $\bd$ and $q_j^{\mp h_j}$, respectively. 
\end{thm}
\noindent
{\it Remark.}\ This theorem is essentially due to Jing\cite{Jing}. However, in \cite{Jing} 
no Serre relations are written explicitly. 
There are also some misprints in Theorem 2.2 in \cite{Jing}:
\begin{itemize}
\item $(1-q^{\mp}w/z)^{-(\al_i|\al_j)/k}_{q^{2k}} $ should be read as $(1-q^{\mp}w/z)^{(\al_i|\al_j)/k}_{q^{2k}}$
\item $(1-q^{\mp}z/w)^{-(\al_i|\al_j)/k}_{q^{2k} }$ should be read as $(1-q^{\mp}z/w)^{(\al_i|\al_j)/k}_{q^{2k}}$
\item $(1-w/z)^{(\al_i|\al_j)/k}_{q^{2k} }$ should be read as $(1-w/z)^{-(\al_i|\al_j)/k}_{q^{2k}}$
\item $(1-z/w)^{(\al_i|\al_j)/k}_{q^{2k}} $ should be read as $(1-z/w)^{-(\al_i|\al_j)/k}_{q^{2k}}$
\end{itemize}

\begin{dfn}\label{z-algebra}
For $k\in \C^\times$ and $(V,\bar{\pi})\in C_k$, we call 
the subalgebra of $\End_\C{V}$
 generated by $Z_{i,m}^\pm(V)$, $q_i^{\pm h_i}\ (i\in I, m\in \Z)$ and $\bd$ 
the quantum $Z$-algebra ${Z}_{V}$ 
 associated with  $(V,\bar{\pi})$. We also define the universal quantum $Z$ algebra 
 $Z_k$ as a topological algebra over $\C[[q^{2k}]]$ in the same way as $\cZ_k$ in 
 Definition \ref{univz-algebra}. We denote the generators in $Z_k$ by $Z_{j,m}^\pm\ (j\in I)$.   
\end{dfn}

\end{appendix}

\renewcommand{\baselinestretch}{0.7}

\end{document}